\documentclass[12pt]{amsart}
\usepackage{mathrsfs}
\usepackage{cases}
\usepackage{epic}
\usepackage{amsfonts}
\usepackage{tikz}
\usepackage{tikz-cd}
\usepackage{graphicx}
\usepackage{amsmath}
\usepackage{amssymb, upgreek}
\usepackage{bm}
\usepackage{latexsym,todonotes}
\usepackage{pdflscape}
\usepackage[all]{xypic}
\usepackage[all]{xy}
\usepackage{tikz-cd}
\usepackage{color}
\usepackage{colordvi}
\usepackage{multicol}
\usepackage{mathtools}
\usepackage[normalem]{ulem}
\usepackage[linktocpage=true]{hyperref}
\hypersetup{colorlinks,linkcolor=blue,urlcolor=cyan,citecolor=blue}
\usepackage{graphicx}
\NewDocumentCommand{\sslash}{s}{%
	\IfBooleanTF{#1}
	{\big/\mkern-7mu\big/}
	{/\mkern-6mu/}%
}
\makeatletter
\newsavebox{\@brx}
\newcommand{\llangle}[1][]{\savebox{\@brx}{\(\m@th{#1\langle}\)}%
	\mathopen{\copy\@brx\kern-0.5\wd\@brx\usebox{\@brx}}}
\newcommand{\rrangle}[1][]{\savebox{\@brx}{\(\m@th{#1\rangle}\)}%
	\mathclose{\copy\@brx\kern-0.5\wd\@brx\usebox{\@brx}}}
\makeatother

\tikzcdset{scale cd/.style={every label/.append style={scale=#1},
		cells={nodes={scale=#1}}}}

\DeclareMathAlphabet{\mathbbb}{U}{bbold}{m}{n}
\DeclareMathOperator{\dimv}{\underline{\dim}}

\makeatletter

\newcommand{\Rmnum}[1]{\textup{\expandafter\@slowromancap\romannumeral #1@}}
\makeatother

\allowdisplaybreaks

\def \tUU{(\tU\otimes\tU)}
\def \tUUi {(\tU\otimes\tU)^\imath}

\begin{document}
	\input xy
	\xyoption{all}
	\newcommand{\iLa}{\Lambda^{\imath}}
	\newcommand{\iadd}{\operatorname{iadd}\nolimits}
	\renewcommand{\mod}{\operatorname{mod}\nolimits}
	\newcommand{\fproj}{\operatorname{f.proj}\nolimits}
	\newcommand{\Fac}{\operatorname{Fac}\nolimits}
	\newcommand{\ci}{{\I}_{\btau}}
	\newcommand{\proj}{\operatorname{proj}\nolimits}
	\newcommand{\inj}{\operatorname{inj}\nolimits}
	\newcommand{\rad}{\operatorname{rad}\nolimits}
	\newcommand{\Span}{\operatorname{Span}\nolimits}
	\newcommand{\soc}{\operatorname{soc}\nolimits}
	\newcommand{\ind}{\operatorname{inj.dim}\nolimits}
	\newcommand{\Ginj}{\operatorname{Ginj}\nolimits}
	\newcommand{\res}{\operatorname{res}\nolimits}
	\newcommand{\np}{\operatorname{np}\nolimits}
	\newcommand{\Mor}{\operatorname{Mor}\nolimits}
	\newcommand{\Mod}{\operatorname{Mod}\nolimits}
	\newcommand{\End}{\operatorname{End}\nolimits}
	\newcommand{\lf}{\operatorname{l.f.}\nolimits}
	\newcommand{\Iso}{\operatorname{Iso}\nolimits}
	\newcommand{\Aut}{\operatorname{Aut}\nolimits}
	\newcommand{\Rep}{\operatorname{Rep}\nolimits}
	
	\newcommand{\colim}{\operatorname{colim}\nolimits}
	\newcommand{\gldim}{\operatorname{gl.dim}\nolimits}
	\newcommand{\cone}{\operatorname{cone}\nolimits}
	\newcommand{\rep}{\operatorname{rep}\nolimits}
	\newcommand{\Ext}{\operatorname{Ext}\nolimits}
	\newcommand{\Tor}{\operatorname{Tor}\nolimits}
	\newcommand{\Hom}{\operatorname{Hom}\nolimits}
	\newcommand{\Top}{\operatorname{top}\nolimits}
	\newcommand{\Coker}{\operatorname{Coker}\nolimits}
	\newcommand{\thick}{\operatorname{thick}\nolimits}
	\newcommand{\rank}{\operatorname{rank}\nolimits}
	\newcommand{\Gproj}{\operatorname{Gproj}\nolimits}
	\newcommand{\Len}{\operatorname{Length}\nolimits}
	\newcommand{\RHom}{\operatorname{RHom}\nolimits}
	\renewcommand{\deg}{\operatorname{deg}\nolimits}
	\renewcommand{\Im}{\operatorname{Im}\nolimits}
	\newcommand{\Ker}{\operatorname{Ker}\nolimits}
	\newcommand{\Coh}{\operatorname{Coh}\nolimits}
	\newcommand{\Id}{\operatorname{Id}\nolimits}
	\newcommand{\Qcoh}{\operatorname{Qch}\nolimits}
	\newcommand{\CM}{\operatorname{CM}\nolimits}
	\newcommand{\sgn}{\operatorname{sgn}\nolimits}
	\newcommand{\utMH}{\operatorname{\cm\ch(\iLa)}\nolimits}
	\newcommand{\GL}{\operatorname{GL}}
	\newcommand{\Perv}{\operatorname{Perv}}
	
	\newcommand{\IC}{\operatorname{IC}}
	\def \hU{\widehat{\U}}
	\def \hUi{\widehat{\U}^\imath}
	\newcommand{\bb}{\psi_*}
	\newcommand{\bvs}{{\boldsymbol{\varsigma}}}
	\def \ba{\mathbf{a}}
	\newcommand{\vs}{\varsigma}
	\def \bfk {\mathbf{k}}

	\def \bd{\mathbf{d}}
	\newcommand{\e}{{\bf 1}}
	\newcommand{\EE}{E^*}
	\newcommand{\dbl}{\operatorname{dbl}\nolimits}
	\newcommand{\ga}{\gamma}
	\newcommand{\tM}{\cm\widetilde{\ch}}
	\newcommand{\la}{\lambda}
	
	\newcommand{\For}{\operatorname{{\bf F}or}\nolimits}
	\newcommand{\coker}{\operatorname{Coker}\nolimits}
	\newcommand{\rankv}{\operatorname{\underline{rank}}\nolimits}
	\newcommand{\diag}{{\operatorname{diag}\nolimits}}
	\newcommand{\swa}{{\operatorname{swap}\nolimits}}
	\newcommand{\supp}{{\operatorname{supp}}}
	
	\renewcommand{\Vec}{{\operatorname{Vec}\nolimits}}
	\newcommand{\pd}{\operatorname{proj.dim}\nolimits}
	\newcommand{\gr}{\operatorname{gr}\nolimits}
	\newcommand{\id}{\operatorname{id}\nolimits}
	\newcommand{\aut}{\operatorname{Aut}\nolimits}
	\newcommand{\Gr}{\operatorname{Gr}\nolimits}
	
	\newcommand{\pdim}{\operatorname{proj.dim}\nolimits}
	\newcommand{\idim}{\operatorname{inj.dim}\nolimits}
	\newcommand{\Gd}{\operatorname{G.dim}\nolimits}
	\newcommand{\Ind}{\operatorname{Ind}\nolimits}
	\newcommand{\add}{\operatorname{add}\nolimits}
	\newcommand{\pr}{\operatorname{pr}\nolimits}
	\newcommand{\oR}{\operatorname{R}\nolimits}
	\newcommand{\oL}{\operatorname{L}\nolimits}
	\def \brW{\mathrm{Br}(W_\btau)}
	\newcommand{\Perf}{{\mathfrak Perf}}
	\newcommand{\cc}{{\mathcal C}}
	\newcommand{\gc}{{\mathcal GC}}
	\newcommand{\ce}{{\mathcal E}}
	\newcommand{\calI}{{\mathcal I}}
	\newcommand{\cs}{{\mathcal S}}
	\newcommand{\cf}{{\mathcal F}}
	\newcommand{\cx}{{\mathcal X}}
	\def \fw{{\mathfrak{w}}}
	
	\newcommand{\cy}{{\mathcal Y}}
	\newcommand{\ct}{{\mathcal T}}
	\newcommand{\cu}{{\mathcal U}}
	\newcommand{\cv}{{\mathcal V}}
	\newcommand{\cn}{{\mathcal N}}
	\newcommand{\mcr}{{\mathcal R}}
	\newcommand{\ch}{{\mathcal H}}
	\newcommand{\ca}{{\mathcal A}}
	\newcommand{\cb}{{\mathcal B}}
	\newcommand{\cj}{{\mathcal J}}
	\newcommand{\cl}{{\mathcal L}}
	\newcommand{\cm}{{\mathcal M}}
	\newcommand{\cp}{{\mathcal P}}
	\newcommand{\cg}{{\mathcal G}}
	\newcommand{\cw}{{\mathcal W}}
	\newcommand{\co}{{\mathcal O}}
	\newcommand{\cq}{{\mathcal Q}}
	\newcommand{\cd}{{\mathcal D}}
	\newcommand{\ck}{{\mathcal K}}
	\newcommand{\calr}{{\mathcal R}}
	\newcommand{\cz}{{\mathcal Z}}
	\newcommand{\ol}{\overline}
	\newcommand{\ul}{\underline}
	\newcommand{\st}{[1]}
	\newcommand{\ow}{\widetilde}
	\renewcommand{\P}{\mathbf{P}}
	\newcommand{\pic}{\operatorname{Pic}\nolimits}
	\newcommand{\Spec}{\operatorname{Spec}\nolimits}
	\newcommand{\Fr}{\mathrm{Fr}}
	\newcommand{\Gp}{\mathrm{Gp}}

	\newtheorem{innercustomthm}{{\bf Theorem}}
	\newenvironment{customthm}[1]
	{\renewcommand\theinnercustomthm{#1}\innercustomthm}
	{\endinnercustomthm}
	
	\newtheorem{innercustomcor}{{\bf Corollary}}
	\newenvironment{customcor}[1]
	{\renewcommand\theinnercustomcor{#1}\innercustomcor}
	{\endinnercustomthm}
	
	\newtheorem{innercustomprop}{{\bf Proposition}}
	\newenvironment{customprop}[1]
	{\renewcommand\theinnercustomprop{#1}\innercustomprop}
	{\endinnercustomthm}
	
	\newtheorem{theorem}{Theorem}[section]
	\newtheorem{acknowledgement}[theorem]{Acknowledgement}
	\newtheorem{algorithm}[theorem]{Algorithm}
	\newtheorem{axiom}[theorem]{Axiom}
	\newtheorem{case}[theorem]{Case}
	\newtheorem{claim}[theorem]{Claim}
	\newtheorem{conclusion}[theorem]{Conclusion}
	\newtheorem{condition}[theorem]{Condition}
	\newtheorem{conjecture}[theorem]{Conjecture}
	\newtheorem{construction}[theorem]{Construction}
	\newtheorem{corollary}[theorem]{Corollary}
	\newtheorem{criterion}[theorem]{Criterion}
	\newtheorem{definition}[theorem]{Definition}
	\newtheorem{example}[theorem]{Example}
	\newtheorem{assumption}[theorem]{Assumption}
	\newtheorem{lemma}[theorem]{Lemma}
	\newtheorem{notation}[theorem]{Notation}
	\newtheorem{problem}[theorem]{Problem}
	\newtheorem{proposition}[theorem]{Proposition}
	\newtheorem{solution}[theorem]{Solution}
	\newtheorem{summary}[theorem]{Summary}
	\newtheorem{hypothesis}[theorem]{Hypothesis}
	\newtheorem*{thm}{Theorem}
	
	\theoremstyle{remark}
	\newtheorem{remark}[theorem]{Remark}
	
	\def \Br{\mathrm{Br}}
	\newcommand{\tK}{K}
	
	\newcommand{\tk}{\widetilde{k}}
	\newcommand{\tU}{\widetilde{{\mathbf U}}}
	\newcommand{\Ui}{{\mathbf U}^\imath}
	\newcommand{\tUi}{\widetilde{{\mathbf U}}^\imath}
	\newcommand{\qbinom}[2]{\begin{bmatrix} #1\\#2 \end{bmatrix} }
	\newcommand{\ov}{\overline}
	\newcommand{\tMHg}{\operatorname{\widetilde{\ch}(Q,\btau)}\nolimits}
	\newcommand{\tMHgop}{\operatorname{\widetilde{\ch}(Q^{op},\btau)}\nolimits}
	
	\newcommand{\rMHg}{\operatorname{\ch_{\rm{red}}(Q,\btau)}\nolimits}
	\newcommand{\dg}{\operatorname{dg}\nolimits}
	\def \fu{{\mathfrak{u}}}
	\def \fv{{\mathfrak{v}}}
	\def \sqq{{\mathbbb{v}}}
	\def \bp{{\mathbf p}}
	\def \bv{{\mathbf v}}
	\def \bw{{\mathbf w}}
	\def \bA{{\mathbf A}}
	\def \bL{{\mathbf L}}
	\def \bF{{\mathbf F}}
	\def \bS{{\mathbf S}}
	\def \bC{{\mathbf C}}
	\def \bU{{\mathbf U}}
	\def \U{{\mathbf U}}
	\def \btau{\varpi}
	\def \La{\Lambda}
	\def \Res{\Delta}
	\newcommand{\ev}{\bar{0}}
	\newcommand{\odd}{\bar{1}}
	\def \fk{\mathfrak{k}}
	\def \ff{\mathfrak{f}}
	\def \fp{{\mathfrak{P}}}
	\def \fg{\mathfrak{g}}
	\def \fn{\mathfrak{n}}
	\def \gr{\mathfrak{gr}}
	\def \Z{\mathbb{Z}}
	\def \F{\mathbb{F}}
	\def \D{\mathbb{D}}
	\def \C{\mathbb{C}}
	\def \N{\mathbb{N}}
	\def \Q{\mathbb{Q}}
	\def \G{\mathbb{G}}
	\def \P{\mathbb{P}}
	\def \K{\mathbb{K}}
	\def \E{\mathbb{K}}
	\def \I{\mathbb{I}}
	
	\def \eps{\varepsilon}
	\def \BH{\mathbb{H}}
	\def \btau{\varrho}
	\def \cv{\varpi}
	
	\def \tR{\widetilde{\bf R}}
	\def \tRZ{\widetilde{\bf R}_\cz}
	\def \hR{\widehat{\bf R}}
	\def \hRZ{\widehat{\bf R}_\cz}
	\def\tRi{\widetilde{\bf R}^\imath}
	\def\hRi{\widehat{\bf R}^\imath}
	\def\tRiZ{\widetilde{\bf R}^\imath_\cz}
	\def\reg{\mathrm{reg}}
	
	\def\hRiZ{\widehat{\bf R}^\imath_\cz}
	\def \tTT{\widetilde{\mathbf{T}}}
	\def \TT{\mathbf{T}}
	\def \br{\mathbf{r}}
	\def \bp{{\mathbf p}}
	\def \tS{\texttt{S}}
	\def \bq{{\bm q}}
	\def \bvt{{v}}
	\def \bs{{ r}}
	\def \tt{{v}}
	\def \k{k}
	\def \bnu{\bm{\nu}}
	\def\bc{\mathbf{c}}
	\def \ts{\textup{\texttt{s}}}
	\def \tt{\textup{\texttt{t}}}
	\def \tr{\textup{\texttt{r}}}
	\def \tc{\textup{\texttt{c}}}
	\def \tg{\textup{\texttt{g}}}
	\def \bW{\mathbf{W}}
	\def \bV{\mathbf{V}}

	\newcommand{\browntext}[1]{\textcolor{brown}{#1}}
	\newcommand{\greentext}[1]{\textcolor{green}{#1}}
	\newcommand{\redtext}[1]{\textcolor{red}{#1}}
	\newcommand{\bluetext}[1]{\textcolor{blue}{#1}}
	\newcommand{\brown}[1]{\browntext{ #1}}
	\newcommand{\green}[1]{\greentext{ #1}}
	\newcommand{\red}[1]{\redtext{ #1}}
	\newcommand{\blue}[1]{\bluetext{ #1}}
	\numberwithin{equation}{section}
	\renewcommand{\theequation}{\thesection.\arabic{equation}}
	
	\newcommand{\wtodo}{\rightarrowdo[inline,color=orange!20, caption={}]}
	\newcommand{\lutodo}{\rightarrowdo[inline,color=green!20, caption={}]}
	\def \tT{\widetilde{\mathcal T}}
	
	\def \tTL{\tT(\iLa)}
	\def \iH{\widetilde{\ch}}
	
	
	\title[Quiver varieties and dual canonical bases]{Quiver varieties and dual canonical bases}

	\author[Ming Lu]{Ming Lu}
	\address{Department of Mathematics, Sichuan University, Chengdu 610064, P.R.China}
	\email{luming@scu.edu.cn}

	\author[Xiaolong Pan]{Xiaolong Pan}
	\address{Department of Mathematics, Sichuan University, Chengdu 610064, P.R.China}
	\email{xiaolong\_pan@stu.scu.edu.cn}


	\maketitle
	\begin{quote}\begin{center}
			{\em Dedicated to Professor Claus Michael Ringel on the occasion of his 80th birthday}
		\end{center}
	\end{quote}
	
	\begin{abstract}
		We survey some recent developments on the theory of dual canonical bases for quantum groups and $\imath$quantum groups. The $\imath$quiver algebras were introduced by Wang and the first author, which are used to give two realizations of quasi-split $\imath$quantum groups of type ADE: one via the $\imath$Hall algebras and the other via the quantum Grothendieck rings of Nakajima-Keller-Scherotzke quiver varieties. The geometric construction of the $\imath$quantum groups produces their dual canonical bases with positivity, generalizing Qin's geometric realization of quantum groups of type ADE. Recently, the authors provided a new construction of the dual canonical basis in the setting of $\imath$Hall algebras, and proved that it is invariant under braid group actions, and obtained the positivity of the transition matrix coefficients from the Hall basis to the dual canonical basis.
		As quantum groups can be regarded as $\imath$quantum groups of diagonal type, we demonstrate that the dual canonical bases of quantum groups coincide with the double canonical bases defined by Berenstein and Greenstein, and resolve several conjectures therein.
	\end{abstract}

	\setcounter{tocdepth}{1}
	
	\tableofcontents
	
	\section{Introduction}




	\subsection{Hall algebras}

	\subsubsection{Quantum groups and Hall algebras}
	
	Ringel \cite{Rin90} introduced a Hall algebra for the representations of a Dynkin quiver over a finite field $\bfk=\F_q$, and proved that this Hall algebra is isomorphic to the positive part $\mathbf{U}^{+}$ of the quantum group $\U$ associated to the underlying Dynkin diagram. Later, Green \cite{Gr95} generalized Ringel's result to arbitrary quantum groups; for a survey see \cite{Sch12}, and we review the basics of Hall algebras in \S\ref{subsec:Hall}. 
	
	It is a natural question to realize the entire quantum groups $\U$ by using Hall algebras. There are many attempts to solve this problem; see e.g. \cite{Kap97,PX97,To06,XX08}. It was eventually solved by 
	Bridgeland \cite{Br13} who constructed a kind of Hall algebras by considering $2$-periodic complexes of quiver representations; cf. Corollary \ref{cor:bridgeland}. In fact, Bridgeland realized the Drinfeld double quantum groups $\tU$, with the Cartan parts doubled. 
	
	Bridgeland's construction has found further generalizations and improvements; see \cite{Gor18} for semi-derived Hall algebras on Frobenius categories, and see \cite{LP21} for semi-derived Ringel-Hall algebras on the categories of complexes over hereditary abelian categories. 
	
	\subsubsection{$\imath$Quantum groups and $\imath$Hall algebras}

	Starting with Satake diagrams, Letzter \cite{Let99,Ko14} introduced quantum symmetric pairs $(\U, \Ui)$, where $\Ui =\Ui_\bvs$ (called an $\imath$quantum group) is a coideal subalgebra of a Drinfeld-Jimbo quantum group $\U$ depending on parameters $\bvs$. The first author and Wang introduced in \cite{LW22a} the universal (quasi-split) $\imath$quantum group  $\tUi$, which is a coideal subalgebra of the Drinfeld double $\tU$, and the $\imath$quantum groups $\Ui_\bvs$ for various parameters can be obtained from $\tUi$ by central reductions. In this survey, we shall restrict ourselves to the universal $\tUi$.  
	
	Drinfeld-Jimbo quantum groups can be viewed as $\imath$quantum groups associated to diagonal Satake diagrams (recall a diagonal Satake diagram is a union of 2 identical copies of a Dynkin diagram with an involution swapping corresponding vertices). The basics on quantum groups and $\imath$quantum groups are reviewed in Section~\ref{sec:QG}.
	
	Building upon the framework of semi-derived Ringel-Hall algebras for $1$-Gorenstein algebras (see \cite{LP21} and \cite[Appendix]{LW22a}), the first author and Wang \cite{LW22a,LW23} formulated the $\imath$Hall algebras of $\imath$quiver algebras associated to $\imath$quivers $(Q,\varrho)$, where $\varrho$ is an involution of the quiver $Q$. The $\imath$Hall algebras of $\imath$quiver algebras are used to provide a realization of the universal quasi-split $\imath$quantum groups $\tU^{\imath}$. For a survey see \cite{LW24}. We review the $\imath$quiver algebras and $\imath$Hall algebras in \S\ref{subsec:i-quiver}--\S\ref{subsec:i-Hall}.

	As an application, BGP type reflection functors are constructed for $\imath$quiver algebras, which induce the relative braid group symmetries of $\tUi$; see \cite{LW21b,LW22b}. The BGP type reflection functors are reviewed in \S\ref{subsec:iHA reflection}.
	
	\subsection{Quiver varieties}

	\subsubsection{Quantum groups and quiver varieties}
	Lusztig in \cite{Lus90,Lus91,Lus93} used perverse sheaves on the varieties of representations of a quiver $Q$ to realize $\U^+$.
	The relation between Ringel’s realization and Lusztig’s categorification is given by sheaf-function correspondence; see \cite{Lus98}. 
	Nakajima \cite{Na01,Na04} further developed Lusztig's construction, and introduced (graded) Nakajima quiver varieties. We review Lusztig's geometric realization in \S\ref{subsec:Lus-variety}.
	
	Via quantum Grothendieck rings of cyclic quiver varieties, Qin \cite{Qin} provided a geometric construction of the Drinfeld doubles of quantum groups of type ADE. Qin’s work was built on the construction of Hernandez--Leclerc \cite{HL15} (who realized
	half a quantum group using graded Nakajima quiver varieties) as well as the concept of
	quantum Grothendieck rings introduced by Nakajima \cite{Na04} and Varagnolo--Vasserot \cite{VV}. These constructions are reviewed in \S\ref{def:RS for QG}.
	
	\subsubsection{$\imath$Quantum groups and quiver varieties}

	As a further generalization of graded Nakajima quiver varieties, motivated by the works \cite{Na01,HL15,LeP13},
	Keller and Scherotzke \cite{KS16,Sch19} formulated the notion of regular/singular Nakajima categories $\mcr, \cs$ for any acyclic quiver $Q$. 
	Note that $\mcr, \cs$ were called the generalized Nakajima categories in \cite{KS16,Sch19}, and they are called Nakajima-Keller-Scherotzke (NKS) categories in \cite{LW21b}. The NKS varieties are by definition (cf. Definition~\ref{def:NKS}) the representation varieties of $\mcr$ and $\cs$. The NKS categories and varieties are reviewed in \S\ref{subsec:NKS-variety}.

	For a Dynkin $\imath$quiver $(Q,\varrho)$, we find that its $\imath$quiver algebra $\Lambda^\imath$ can be realized  as a singular NKS category $\cs^\imath$, and then use the (dual) quantum Grothendieck ring of perverse sheaves over $\cs^\imath$ to realize quasi-split $\imath$quantum group $\tUi$; see \cite{LW21b}. Qin's construction \cite{Qin} of the Drinfeld double $\tU$ uses such cyclic quiver varieties, which can be viewed as NKS varieties of $\cs^\imath$ of diagonal type; cf. the $\imath$quiver algebras are used to formulate Bridgeland's Hall algebras; cf. \cite{LW22a}. The geometric realization of $\imath$quantum groups is reviewed in \S\ref{subsec:geometric-iQG}.

	\subsection{Dual canonical bases}

	\subsubsection{Canonical bases of quantum groups}
	
	Based on the geometric realization of $\U^+$ in \cite{Lus90}, Lusztig produced the canonical basis of $\U^+$ by simple perverse sheaves; cf. also \cite{Ka91} for the crystal basis. Lusztig \cite{Lus90} also introduced braid group symmetries on the quantum group $\U$ and used them to define the PBW basis. He also used the PBW basis to give an elementary construction of the canonical basis of $\U^+$; also see \cite[\S11.6]{DDPW} by noting that PBW bases in this case coincide with Hall bases (i.e. natural bases of Hall algebras). The canonical basis of $\U^+$ via perverse sheaves is reviewed in \S\ref{subsec:Lus-variety}, and via Hall basis is reviewed in \S\ref{subsec: CB by HA}. 
	
	The dual canonical basis of $\U^+$ also has very nice properties and has deep connections to cluster algebras. We used Hall basis to give it an elementary construction; see \cite{LP25}. This is reviewed in \S\ref{dCB of HA subsec}.

	\subsubsection{Dual canonical bases via perverse sheaves}

	It is a natural question to expand the canonical basis on $\U^+$ to the whole quantum group $\U$. Lusztig found that this is impossible, and he found a way to solve this question by introducing a modified version $\dot{\U}$ of the quantum group; see \cite{Lus93}. Inspired by Lusztig's construction, Bao and Wang \cite{BW18,BW18b} found the canonical bases on $\dot{\U}^\imath$, the modified $\imath$quantum group.
	
	A bonus of the
	geometric approach \cite{Qin,LW21b} is the construction of an integral and positive basis on $\tU$ and $\tUi$, which contains as subsets
	the (mildly rescaled) dual canonical bases of Lusztig for $\U^+$
	\cite{Lus90}. For the quantum group, this positive basis expands the dual canonical bases of $\U^\pm$ to the whole quantum group, and solved the natural question above in some sense. We call this positive basis the dual canonical basis of $\tU$ and $\tUi$, which corresponds to the perverse sheaves on the NKS quiver varieties. 

	\subsubsection{Dual canonical bases via $\imath$Hall algebras}
	Recently, we \cite{LP25,LP26a,LP26b} used $\imath$Hall algebras to reconstruct the dual canonical bases on quantum groups and $\imath$quantum groups. This is based on an anti-involution $\ov{\cdot}$ of $\tUi$ (called bar involution) and Lusztig's Lemma \cite{BZ14}; note that it is different to Lusztig's bar involution \cite{Lus90}. This is reviewed in \S\ref{subsec:dCB-iHall}.
	
	In \cite{LP25}, we also constructed Fourier transforms on $\imath$Hall algebras for different orientations of the quiver $Q$, and proved that the Fourier transforms preserve the dual canonical basis. By using the BGP type reflection functors of $\imath$quiver algebras, we proved the dual canonical bases are invariant (as a set) under the relative braid symmetries. These constructions are reviewed briefly in \S\ref{subsec:iHA reflection}.
	
	We also computed the dual canonical bases for the quantum group and $\imath$quantum group of rank $1$: $\tU_v(\mathfrak{sl}_2)$, $\tUi_v(\mathfrak{sl}_2)$ in \cite{LP26a,LP26b}; see \S\ref{sec:example} for a review.
	
	\subsubsection{Double canonical bases of quantum groups}
	Berenstein and Greenstein \cite{BG17} have constructed a double canonical basis of $\tU$ by combining the dual canonical bases of the positive and negative halves of $\tU$. This basis has many properties similar to the dual canonical basis, for example they are both bar-invariant and contain the dual canonical bases of the positive and negative halves. We \cite{LP26b} obtained the coincidence of the dual canonical basis and the double canonical basis for $\tU$, and then solved several conjectures for double canonical bases therein. We reviewed the double/dual canonical bases of quantum groups in \S\ref{subsec:dCB-QG}.
	
	\subsubsection{Open problems}

	In the final Section~\ref{sec:open}, we formulate several open problems in dual canonical bases and put them in context.

	\vspace{2mm}
	\noindent {\bf Acknowledgement.}
	We thank Shiquan Ruan and Weiqiang Wang for their collaboration on related projects, as well as for their insightful comments and stimulating discussions. We also gratefully acknowledge the anonymous referee for their careful reading and thoughtful corrections.
	This work was partially supported by the National Natural Science Foundation of China (No. 12171333, 12261131498).

	\section{Quantum groups and $\imath$quantum groups}
	\label{sec:QG}
	
	Let $\I=\{1,\dots,n\}$ be the index set. 
	Let $C=(c_{ij})_{i,j \in \I}$ be the Cartan matrix of a simply-laced semi-simple Lie algebra $\fg$.
	Let $\Delta^+=\{\alpha_i\mid i\in\I\}$ be the set of simple roots of $\fg$, and denote the root lattice by $\Z^{\I}:=\Z\alpha_1\oplus\cdots\oplus\Z\alpha_n$. Let $\Phi^+$ be the set of positive roots. The simple reflection $s_i:\Z^{\I}\rightarrow\Z^{\I}$ is defined to be $s_i(\alpha_j)=\alpha_j-c_{ij}\alpha_i$, for $i,j\in \I$.
	Denote the Weyl group by $W =\langle s_i\mid i\in \I\rangle$.

	Let $v$ be an indeterminate. Denote, for $r,m \in \N$,
	\[
	[r]=\frac{v^r-v^{-r}}{v-v^{-1}},
	\quad
	[r]!=\prod_{i=1}^r [i], \quad \qbinom{m}{r} =\frac{[m][m-1]\ldots [m-r+1]}{[r]!}.
	\]

	\subsection{Quantum groups}\label{subsec:QG}
	Following \cite{BG17}, the Drinfeld double $\hU := \hU_v(\fg)$ is defined to be the $\Q(v^{1/2})$-algebra generated by $E_i,F_i, \tK_i,\tK_i'$, $i\in \I$, 
	subject to the following relations:  for $i, j \in \I$,
	\begin{align}
		[E_i,F_j]= \delta_{ij}(v^{-1}-v) (\tK_i-\tK_i'),  &\qquad [\tK_i,\tK_j]=[\tK_i,\tK_j']  =[\tK_i',\tK_j']=0,
		\label{eq:KK}
		\\
		\tK_i E_j=v^{c_{ij}} E_j \tK_i, & \qquad \tK_i F_j=v^{-c_{ij}} F_j \tK_i,
		\label{eq:EK}
		\\
		\tK_i' E_j=v^{-c_{ij}} E_j \tK_i', & \qquad \tK_i' F_j=v^{c_{ij}} F_j \tK_i',
		\label{eq:K2}
	\end{align}
	and for $i\neq j \in \I$,
	\begin{align}
		& \sum_{r=0}^{1-c_{ij}} (-1)^r \left[ \begin{array}{c} 1-c_{ij} \\r \end{array} \right]  E_i^r E_j  E_i^{1-c_{ij}-r}=0,
		\label{eq:serre1} \\
		& \sum_{r=0}^{1-c_{ij}} (-1)^r \left[ \begin{array}{c} 1-c_{ij} \\r \end{array} \right]  F_i^r F_j  F_i^{1-c_{ij}-r}=0.
		\label{eq:serre2}
	\end{align}
	
	We define $\tU=\tU_v(\fg)$ to be the $\Q(v^{1/2})$-algebra constructed from $\widehat{\U}$ by making $\tK_i,\tK_i'$ ($i\in\I$) invertible. Then $\tU$ and $\hU$ are $\Z^\I$-graded by setting $\deg E_i=\alpha_i$, $\deg F_i=-\alpha_i$, $\deg K_i=0=\deg K_i'$. 
	
	\begin{remark}
		Note that $\tU$ defined here is different but isomorphic to the usual one given in \cite{Lus93}. 
		In fact, denote 
		\begin{align}
			\label{eq:Udj-gen}
			\ce_i=\frac{E_i}{v^{-1}-v},\qquad \cf_i=\frac{F_i}{v-v^{-1}},\qquad \forall i\in\I.
		\end{align} 
		Then the presentation of $\tU$ given in \cite{Lus93} is generated by $\ce_i,\cf_i,K_i,K_i'$ ($i\in\I$).
	\end{remark}

	The quantum group $\bU$ is defined to be quotient algebra of $\hU$ (also $\tU$) modulo the ideal generated by $K_iK_i'-1$ ($i\in\I$). In other words, $\U$ is the $\Q(v^{1/2})$-algebra generated by $E_i,F_i, K_i, K_i^{-1}$, $i\in \I$, subject to the  relations modified from \eqref{eq:KK}--\eqref{eq:serre2} with $\tK_i'$ replaced by $K_i^{-1}$.
	
	Let $\hU^+$ be the subalgebra of $\hU$ generated by $E_i$ $(i\in \I)$, $\hU^0$ be the subalgebra of $\widehat{\bU}$ generated by $\tK_i, \tK_i'$ $(i\in \I)$, and $\hU^-$ be the subalgebra of $\widehat{\bU}$ generated by $F_i$ $(i\in \I)$, respectively.
	The subalgebras $\tU^+$, $\tU^0$ and $\tU^-$ of $\tU$, and also $\bU^+$, $\bU^0$ and $\bU^-$ of $\bU$ are defined similarly. Then we have the triangular decompositions:
	\begin{align*}
		\hU=\hU^+\otimes\hU^0\otimes \hU^-,\qquad 
		\widetilde{\bU} =\widetilde{\bU}^+\otimes \widetilde{\bU}^0\otimes\widetilde{\bU}^-,
		\qquad
		\bU &=\bU^+\otimes \bU^0\otimes\bU^-.
	\end{align*}
	Clearly, ${\bU}^\pm\cong\widetilde{\bU}^\pm\cong\hU^\pm$. 
	For any $\mu=\sum_{i\in\I}m_i\alpha_i\in\Z^\I$, we denote by $K_\mu=\prod_{i\in\I} K_i^{m_i}$, $K_\mu'=\prod_{i\in\I} (K_i')^{m_i}$ in $\tU$ (or $\U$); we can view $\tK_\mu,\tK_\mu'$ in $\hU$ if $\mu\in\N^\I$.
	
	\begin{lemma}[cf. \cite{Lus90,BG17}]\label{QG bar-involution def}
		\begin{enumerate}
			\item  There exists an anti-involution $u\mapsto \ov{u}$ on $\hU$ (also $\tU$, $\U$) given by $\ov{v^{1/2}}=v^{-1/2}$, $\ov{E_i}=E_i$, $\ov{F_i}=F_i$, and $\ov{K_i}=K_i$, $\ov{K_i'}=K_i'$, for $i\in\I$.
			\item There exists an involution $u\mapsto \psi(u)$ on $\hU$ (also $\tU$, $\U$) given by $\psi(v^{1/2})=v^{-1/2}$, $\psi(E_i)=E_i$, $\psi(F_i)=F_i$, and $\psi(K_i)=K_i'$, $\psi(K_i')=K_i$, for $i\in\I$.
		\end{enumerate}
	\end{lemma}
	
	The algebras $\hU$ (and $\tU$, $\U$) are Hopf algebras, with the coproduct $\Delta$ defined by
	\begin{align}
		\begin{split}
			\Delta(E_i)=E_i\otimes 1+K_i\otimes E_i,\quad &\Delta(F_i)=1\otimes F_i+F_i\otimes K_i',
			\\
			\Delta(K_i)=K_i\otimes K_i,\quad &\Delta(K_i')=K_i'\otimes K_i'.
		\end{split}
	\end{align}
	
	Let $\Br(W)$ be the braid group associated to the Weyl group $W$, generated by $t_i$ ($i\in\I$).
	Lusztig introduced braid group symmetries on $\U$ \cite[\S37.1.3]{Lus93}. Lusztig's braid group symmetries can be lifted to the Drinfeld double $\tU$; see \cite[Propositions 6.20–6.21]{LW21a}, which are denoted by $\widetilde{T}_i$ (compare with \cite[Theorem 1.13]{BG17}).
	\begin{proposition}
		\label{prop:braid1}
		There exists an automorphism $\widetilde{T}_{i}$, for $i\in \I$, on $\tU$ such that
		\begin{align*}
			&\widetilde{T}_{i}(K_\mu)= K_{s_i(\mu)},
			\qquad \widetilde{T}_{i}(K'_\mu)= K'_{s_i(\mu)},\;\;\forall \mu\in \Z^\I,\\
			&\widetilde{T}_{i}(E_i)=v(K_i')^{-1}F_i,\qquad \widetilde{T}_{i}(F_i)=v^{-1}E_iK_i^{-1},\\
			&\widetilde{T}_i(E_j)=E_j,\qquad \widetilde{T}_i(F_j)=F_j,\qquad \text{ if }c_{ij}=0,
			\\
			&\widetilde{T}_{i}(E_j)=\frac{v^{\frac{1}{2}}E_iE_j-v^{-\frac{1}{2}}E_jE_i}{v-v^{-1}},\qquad  \widetilde{T}_{i}(F_j)=\frac{v^{\frac{1}{2}}F_iF_j-v^{-\frac{1}{2}}F_jF_i}{v-v^{-1}}, \text{ if } c_{ij}=-1. 
		\end{align*}
		Moreover, there exists a group homomorphism $\Br(W)\rightarrow \Aut(\tU)$, $t_i\mapsto \widetilde{T}_i$ for $i\in\I$.    
	\end{proposition}
	
	We can therefore define
	\begin{align}\widetilde{T}_w 
		:= \widetilde{T}_{i_1}\cdots
		\widetilde{T}_{i_r} \in \Aut(\tU),
	\end{align}
	where $w = s_{i_1}\cdots s_{i_r}$ is any reduced expression of $w \in W$. 
	
	\begin{lemma}
		\label{lem:QGbraid-bar}
		The braid group actions $\widetilde{T}_{i}$ commute with the bar-involution, i.e., $\ov{\widetilde{T}_i(u)}=\widetilde{T}_i(\ov{u})$ for any $u\in\tU$.
	\end{lemma}

	\subsection{The $\imath$quantum groups}
	
	For a  Cartan matrix $C=(c_{ij})$, let $\Aut(C)$ be the group of all permutations $\btau$ of the set $\I$ such that $c_{ij}=c_{\btau i,\btau j}$. An element $\btau\in\Aut(C)$ is called an \emph{involution} if $\btau^2=\Id$.

	In this paper, we always assume that $\btau\in\Aut(C)$ is an involution such that $c_{i,\btau i}=0$ for all $i\neq \btau i$. We denote by $\bs_{i}$ the following element of order $2$ in the Weyl group $W$, i.e.,
	\begin{align}
		\label{def:simple reflection}
		\bs_i= \left\{
		\begin{array}{ll}
			s_{i}, & \text{ if } \btau i=i;
			\\
			s_is_{\btau i}, & \text{ if } \btau i\neq i.
		\end{array}
		\right.
	\end{align}
	It is well known (cf., e.g., \cite{KP11}) that the {\rm restricted Weyl group} associated to the quasi-split symmetric pair $(\fg, \fg^\theta)$ can be identified with the following subgroup $W_\btau$ of $W$:
	\begin{align}
		\label{eq:Wtau}
		W_{\btau} =\{w\in W\mid \btau w =w \btau\}
	\end{align}
	where $\btau$ is regarded as an automorphism of the root lattice $\Z^\I$. Moreover, $W_{\btau}$ can be identified with a Weyl group with $\bs_i$ ($i\in \I_\btau$) as its simple reflections.
	
	For a  Cartan matrix $C=(c_{ij})$, 
	let $\btau$ be an involution in $\Aut(C)$. We define the universal $\imath$quantum group ${\hU}^\imath:=\hU'_v(\fk)$ (resp. $\tUi:=\tU'_v(\fk)$) to be the $\Q(v^{1/2})$-subalgebra of $\hU$ (resp. $\tU$) generated by
	\begin{equation}
		\label{eq:Bi}
		B_i= F_i +  E_{\btau i} \tK_i',
		\qquad \tk_i = \tK_i \tK_{\btau i}', \quad \forall i \in \I,
	\end{equation}
	(with $\tk_i$ invertible in $\tUi$).
	Let $\hU^{\imath 0}$ be the $\Q(v^{1/2})$-subalgebra of $\hUi$ generated by $\tk_i$, for $i\in \I$. Similarly, let $\tU^{\imath 0}$ be the $\Q(v^{1/2})$-subalgebra of $\tUi$ generated by $\tk_i^{\pm1}$, for $i\in \I$. 

	A presentation for $\tUi$ can be found in \cite[Proposition~ 6.4]{LW22a}. The following is a modified version based on the quantum groups given in \S\ref{subsec:QG} (this kind of presentation also holds for $\hUi$ by omitting that $\tk_i$ ($i\in\I$) are invertible).
	\begin{proposition}
		\label{prop:Serre}
		The $\Q(v^{1/2})$-algebra $\tUi$ has a presentation with generators $B_i, \tk_i$ $(i\in \I)$, where $\tk_i$ are invertible, subject to the relations \eqref{relation1}--\eqref{relation2}: for $\ell \in \I$, and $i\neq j \in \I$,
		\begin{align}
			&\tk_i \tk_\ell =\tk_\ell \tk_i,
			\quad
			\tk_\ell B_i = v^{c_{\btau \ell,i} -c_{\ell i}} B_i \tk_\ell,
			\label{relation1}
			\\
			&B_iB_{j}-B_jB_i =0, \quad \text{ if } c_{ij} =0 \text{ and }\btau i\neq j,
			\label{relationBB}
			\\
			&\sum_{s=0}^{1-c_{ij}} (-1)^s \qbinom{1-c_{ij}}{s} B_i^{s}B_jB_i^{1-c_{ij}-s} =0, \quad \text{ if } j \neq \btau i\neq i,
			\\
			&B_{\btau i}B_i -B_i B_{\btau i} =  (v^{-1}-v) (\tk_i -\tk_{\btau i}),
			\quad \text{ if } \btau i \neq i,
			\label{relation5}
			\\
			&B_i^2B_{j} - [2] B_iB_{j}B_i +B_{j}B_i^2 = -v(v-v^{-1})^2 \tk_i B_{j},
			%
			\quad \text{ if }  c_{ij} = -1,\btau i=i.
			\label{relation2}
		\end{align}
	\end{proposition}
	
	It is known \cite{Let99, Ko14,LW22a} that the algebra $\widetilde{\bU}^\imath$ (resp. $\hUi$) is a right coideal subalgebra of $\widetilde{\bU}$ (resp. $\hU$); we call $(\widetilde{\bU}, \widetilde{\bU}^\imath)$ and $(\hU,\hUi)$ quantum symmetric pairs. We shall refer to $\hUi$ and $\tUi$ 
	as the universal {\em (quasi-split) $\imath${}quantum groups} (cf. the $\imath$quantum groups defined in  \cite{Let99,Ko14}); they are called {\em split} if $\btau =\Id$. 
	
	From the presentation given in Proposition~\ref{prop:Serre} we can define an anti-involution as follows.
	\begin{lemma}
		\label{iQG bar-involution def}
		There exists an anti-involution $u\mapsto \ov{u}$ on $\hUi$ (also $\tUi$) given by $\ov{v^{1/2}}=v^{-1/2}$, $\ov{B_i}=B_i$, $\ov{\K_i}=\K_i$, for $i\in\I$. In particular, $\ov{\tk_i}=\tk_i$ if $\varrho i\neq i$; $\ov{\tk_i}=v^2\tk_i$ if $\varrho i=i$.
	\end{lemma}
	
	\begin{example}
		\label{ex:QGvsiQG}
		Let us explain that quantum groups as $\imath$quantum groups of diagonal type. 
		Consider the $\Q(v^{1/2})$-subalgebra $\tUUi$ of $\tUU$
		generated by
		\[
		\mathcal{K}_i:=\tK_{i} \tK_{i^{\diamond}}', \quad
		\mathcal{K}_i':=\tK_{i^{\diamond}} \tK_{i}',  \quad
		\cb_{i}:= F_{i}+ E_{i^{\diamond}} \tK_{i}', \quad
		\cb_{i^{\diamond}}:=F_{i^{\diamond}}+ E_{i} \tK_{i^{\diamond}}',
		\qquad \forall i\in \I.
		\]
		Here we drop the tensor product notation and use instead $i^\diamond$ to index the generators of the second copy of $\tU$ in $\tUU$. There exists a $\Q(v^{1/2})$-algebra isomorphism $\widetilde{\phi}: \tU \rightarrow \tUUi$ such that
		\[
		\widetilde{\phi}(E_i)= \cb_{i},\quad \widetilde{\phi}(F_i)= \cb_{i^{\diamond}}, \quad \widetilde{\phi}(\tK_i)= \mathcal{K}_i', \quad \widetilde{\phi}(\tK_i')= \mathcal{K}_i, \qquad \forall  i\in \I.
		\]
	\end{example}

	Choose one representative for each $\btau$-orbit on $\I$, and let
	\begin{align}\label{eq:ci}
		\ci = \{ \text{the chosen representatives of $\btau$-orbits in $\I$} \}.
	\end{align} 
	The braid group associated to the relative Weyl group for $(\fg, \fg^\theta)$ is of the form
	\begin{equation}
		\label{eq:braidCox}
		\brW =\langle \br_i \mid i\in \I_\btau \rangle
	\end{equation}
	where $\br_i$ satisfy the same braid relations as for $\bs_i$ in $W_{\varrho}$ (but no quadratic relations on $r_i$ are imposed). The following result was proved in \cite{LW21a} using $\imath$Hall algebra technique (cf. \cite{KP11, WZ23}), 
	where for any $i\in\I$, we set
	\begin{align}
		\label{eq:bbKi}
		\K_i:=v\tk_i, \text{ if }\varrho i=i;
		\qquad
		\K_j:=\tk_j, \text{ otherwise.}
	\end{align}
	and $\K_\alpha:=\prod_{i\in\I}\K_i^{a_i}$ if $\alpha\in\Z^\I$. For elements $x,y\in\Ui$, we also define $[x,y]_v:=xy-vyx$.
	
	\begin{theorem} [\text{\cite[Theorem 6.8]{LW21a}}]
		\label{thm:Ti}
		\begin{enumerate}
			\item[(1)] For $i\in \I$ such that $\varrho i=i$, there exists an automorphism $\tTT_i$ of the $\Q(v^{1/2})$-algebra $\tUi$ such that
			$\tTT_i(\K_\mu) =\K_{\bs_i\mu}$, and
			\[
			\tTT_i(B_j)= \begin{cases}
				\K_i^{-1} B_i,  &\text{ if }j=i,\\
				B_j,  &\text{ if } c_{ij}=0, \\
				\frac{v^{\frac{1}{2}}B_iB_j-v^{-\frac{1}{2}}B_jB_i}{v-v^{-1}},  & \text{ if }c_{ij}=-1
			\end{cases}
			\]
			for $\mu\in \Z^\I$ and $j\in \I$.
			
			\item[(2)] For $i\in \I$ such that $c_{i,\varrho i}=0$, there exists an automorphism $\tTT_i$ of the $\Q(v^{1/2})$-algebra $\tUi$ such that
			$\tTT_i(\K_\mu)=\K_{\bs_i\mu}$, and
			\[
			\tTT_i(B_j)= \begin{cases}
				\frac{v^{\frac{1}{2}}B_iB_j-v^{-\frac{1}{2}}B_jB_i}{v-v^{-1}},  & \text{ if }c_{ij}=-1 \text{ and } c_{\varrho i,j}=0,
				\\
				\frac{v^{\frac{1}{2}}B_{\varrho i}B_j-v^{-\frac{1}{2}}B_jB_{\varrho i}}{v-v^{-1}},  & \text{ if } c_{ij}=0 \text{ and }c_{\varrho i,j}=-1 ,
				\\
				\frac{v^{-1}\big[[B_j,B_i]_v,B_{\varrho i}\big]_v}{(v-v^{-1})^2}+B_j\K_i  & \text{ if } c_{ij}=-1 \text{ and }c_{\varrho i,j}=-1 ,
				\\
				v\K_{i}^{-1} B_{\btau i},  & \text{ if }j=i,
				\\
				v\K_{\btau i}^{-1} B_i,  &\text{ if }j=\varrho i,
				\\
				B_j, & \text{ otherwise;}
			\end{cases}
			\]
			for $\mu\in \Z^\I$ and $j\in \I$. 
		\end{enumerate}
		Moreover, there exists a homomorphism $\brW \rightarrow \Aut( \tUi)$, $\br_i\mapsto \tTT_i$, for all $i\in \I_\btau$.
	\end{theorem}
	
	Similar to Lemma \ref{lem:QGbraid-bar}, one checks that the $\widetilde{\TT}_{i}$ are equivariant under bar-involution.
	
	\begin{lemma}
		The braid group actions $\widetilde{\TT}_{i}$ commute with the bar-involution, i.e., $\ov{\widetilde{\TT}_i(u)}=\widetilde{\TT}_i(\ov{u})$ for any $u\in\tUi$.
	\end{lemma}

	\section{Hall algebras and Lusztig's quiver varieties}
	
	Let $\bfk =\F_q$ be a finite field of $q$ elements. Let $\mathbbb{v}=\sqrt{q}$. 
	Let $Q=(Q_0=\I,Q_1)$ be a Dynkin quiver, and $\mod(\bfk Q)$ the category of finite-dimensional representations of $Q$ over $\bfk$. We denote by $\dimv M$ the dimension vector of $M\in\mod(\bfk Q)$. 
	Denote by $S_i$ the simple module supported at the vertex for each $i\in\I$. Identify $\dimv S_i$ with the simple root $\alpha_i$. Let $\langle\cdot,\cdot\rangle_Q$ be the Euler form of $\mod(\bfk Q)$, that is,
	$\langle M,N \rangle_Q=\dim\Hom_{\bfk Q}(M,N)-\dim\Ext^1_{\bfk Q}(M,N)$. 
	
	For a finite set $\cs$, we denote by $|\cs|$ its cardinality. 
	
	\subsection{Hall algebras}
	\label{subsec:Hall}
	
	Denote by $\Iso(\mod(\bfk Q))$ the set of isoclasses of $\mod(\bfk Q)$. We use $[M]$ to denote the isoclass of $M$. For any $X,Y,Z\in\mod(\bfk Q)$, we set $\Ext^1_{\bfk Q}(X,Z)_Y\subseteq \Ext^1_{\bfk Q}(X,Z)$ to be the subset parameterizing extensions with the middle term isomorphic to $Y$.
	Let $\widetilde{\ch}(\bfk Q)$ be the Ringel-Hall algebra \cite{Rin90} of $\mod(\bfk Q)$. That is, $\widetilde{\ch}(\bfk Q)$ is defined on the $\Q(\sqq)$-vector space whose basis is formed by the isomorphism classes $[X]$ of objects $X$ of $\mod(\bfk Q)$, with the multiplication
	defined by
	\begin{align}
		\label{eq:mult}
		[X]\cdot [Z]=\sqq^{\langle X,Z\rangle_Q}\sum_{[Y]\in \Iso(\mod(\bfk Q))}G_{XZ}^Y[Y],\text{ where }G_{XZ}^Y:=\frac{|\Ext_{\bfk Q}^1(X,Z)_Y|}{|\Hom_{\bfk Q}(X,Z)|}.
	\end{align}
	It is well known that the algebra $\widetilde{\ch}(\bfk Q)$ is associative and unital; see \cite{Rin90,Br13}.
	
	We denote by $[\![X]\!]:=\frac{[X]}{a_X}$ for any $X\in\mod(\bfk Q)$. For any three objects $X,Y,Z$, let
	\begin{align}
		\label{eq:Fxyz}
		F_{XZ}^Y&:= \big |\{L\subseteq Y \mid L \cong Z,  Y/L\cong X\} \big |.
	\end{align}
	Let $\Aut(X)$ be the automorphism group of $X\in\mod(\bfk Q)$, and denote by $a_X=|\Aut(X)|$.  The Riedtmann-Peng formula states that
	\[
	F_{XZ}^Y= G_{XZ}^Y \cdot \frac{a_Y}{a_X a_Z}
	\] 
	so we have
	\begin{align}
		[\![X]\!]\cdot[\![Z]\!]=\sum_{[\![Y]\!]}F_{XZ}^Y[\![Y]\!].
	\end{align}
	
	For any $\beta\in\Phi^+$, by Gabriel Theorem, there exists a unique indecomposable $\bfk Q$-module $M_q(\beta)$ (up to isomorphism)  such that $\dimv M_q(\beta)=\beta$. 
	Recall that $\mathfrak{P}:=\mathfrak{P}(Q)$ is the set of functions $\lambda: \Phi^+\rightarrow \N$. We denote by $$M_q(\lambda)=\bigoplus_{\beta\in\Phi^+}M_q(\beta)^{\oplus\lambda(\beta)},\qquad \forall \lambda\in\mathfrak{P}.$$
	
	From \cite{Rin90}, we know for any $\lambda,\mu,\nu\in\mathfrak{P}$, there exists a polynomial $f^{\lambda}_{\mu,\nu}(v)\in\Z[v,v^{-1}]$ such that
	$$f^{\lambda}_{\mu,\nu}(\sqq)=\sqq^{\langle M_q(\mu),M_q(\nu)\rangle_Q}F_{M_q(\mu),M_q(\nu)}^{M_q(\lambda)}.$$
	
	We let $\widetilde{\ch}(Q)$ be the generic Hall algebra of $Q$. It is a $\Q(v)$-vector space with basis $\{\fw_\lambda\mid\lambda\in\mathfrak{P}\}$, and the multiplication is defined by 
	\[\fw_\mu\cdot\fw_\nu=\sum_{\lambda}f^{\lambda}_{\mu,\nu}(v)\fw_\lambda.\]
	
	Let $C=(c_{ij})_{i,j\in\I}$ be the Cartan matrix of $Q$, that is, $c_{ij}=\langle S_i,S_j\rangle_Q+\langle S_j,S_i\rangle_Q$, and $\tU$ be the associated quantum group (defined over $\Q(v)$). Then the following theorem is well known.

	\begin{theorem}[\cite{Rin90}]\label{thm: iso R^pm}
		For a Dynkin quiver $Q$, there exists isomorphisms of $\Q(v)$-algebras:
		\begin{align*}
			R^+:\U^+\longrightarrow \widetilde{\ch}(Q),\quad \ce_i\mapsto \fw_{\alpha_i},\quad\forall i\in\I;
			\\
			R^-:\U^-\longrightarrow \widetilde{\ch}(Q),\quad \cf_i\mapsto \fw_{\alpha_i},\quad\forall i\in\I.
		\end{align*}    
	\end{theorem}

	\subsection{Lusztig's quiver varieties}
	\label{subsec:Lus-variety}
	
	Let $Q=(Q_0,Q_1)$ be a quiver and set $\I=Q_0$, $\Omega=Q_1$. Denote by $s,t:Q_1\rightarrow Q_0$ the source and target maps, respectively. For $\bd\in\N^\I$ and an $\I$-graded space $V=\oplus_{i\in\I}V_i$ of dimension vector $\bd$, we define the representation space to be
	\begin{equation}\label{def:Lusztig space E_V}
		E_V=\bigoplus_{h \in \Omega}\Hom(V_{s(h)},V_{t(h)}).
	\end{equation}
	Any point $x\in E_V$ is a representation of $Q$. 
	Note that the group $G_V:=\prod_i\GL(V_i)$ acts on $E_V$ by 
	\[(g\cdot x)_h=g_{t(h)} x_h g_{s(h)}^{-1}\]
	and its orbits parametrize isoclasses of $\mod(\bfk Q)$ of dimension vector $\bd$. For any $x\in E_V$, we denote by $\mathfrak{O}_x$ its orbit. 
	
	Let $\bd\in\N^\I$ and fix an $\I$-graded vector space $V$ of dimension vector $\bd$. We define $I_\bd$ to be the set of sequences $\nu=(\nu^1,\dots,\nu^m)$ satisfying $\bd=\sum_l\nu^l$ and each $\nu^l=n_l\alpha_l$ for some integer $n_l>0$ and simple root $\alpha_l$. For $\nu\in I_\bd$, a \textit{flag of type $\nu$ in $V$} is defined to be a sequence of $\I$-graded subspaces
	\[V=V^0 \supset V^1 \supset \cdots \supset V^n=\{0\}\]
	such that each $V^{l-1}/V^{l}$ has dimension $\nu^l$. Let $\mathcal{F}_{\nu}$ be the variety given by
	\[\mathcal{F}_{\nu}=\{(x,f) \mid \text{$x\in E_V$ and $f$ is an $x$-stable flag of type $\nu$}\}.\]
	Clearly the group $G_V$ acts on $\mathcal{F}_{\nu}$ by $g\cdot (x,f)=(g\cdot x,g(f))$.
	
	\begin{lemma}[\cite{Lus93}]
		The variety $\mathcal{F}_{\nu}$ is smooth and the natural projection $\pi_\nu:\mathcal{F}_{\nu}\to E_V$ is proper and $G_V$-equivariant.    
	\end{lemma}
	
	Let $\cd(E_V):=\cd^b_c(E_V)$ be the bounded derived category of $\bar{\Q}_\ell$-constructible sheaves on $E_V$. Here, $\ell$ is a fixed prime number, and $\bar{\Q}_\ell$ is an algebraic closure of the field of $\ell$-adic numbers. By the decomposition theorem of \cite{BBD}, the complex
	\[L_{\nu} = (\pi_{\nu})_! (\overline{\Q}_\ell [d(\nu)])\]
	is semisimple, where $d(\nu)$ denotes the dimension of $\mathcal{F}_{\nu}$. Let $\mathcal{P}_V$ be the subcategory of $\cd(E_V)$ consisting of simple perverse sheaves $L$ such that $L[n]$ appears in the direct sum decomposition of some $L_{\nu}$, with $\nu\in I_\bd$ and $n\in\Z$. The category $\mathcal{Q}_V$ is then the subcategory of $\cd(E_V)$ generated by direct sums of shift of complexes in $\mathcal{P}_V$. We refer to complexes in $\mathcal{P}_V$ as \textit{Lusztig's perverse sheaves}.

	
	
	Let $\bd,\bd',\bd''\in\N^\I$ be such that $\bd=\bd'+\bd''$. We fix $\I$-graded vector spaces $V,V',V''$ of dimension vector $\bd,\bd',\bd''$ respectively. 
	
	For any $x\in E_V$, an $\I$-graded vector subspace $W\subset V$ is called $x$-stable if $x_h(W_{s(h)})\subset W_{t(h)}$ for any $h\in\Omega$. Let $E_2$ be the variety of $(x,W)$, where $x\in E_V$ and $W\subset V$ is an $x$-stable, $\I$-graded vector subspace of dimension vector $\bd''$; and $E_1$ be the variety of $(x,W,\rho_1,\rho_2)$, where $(x,W)\in E_2$ and $\rho':V/W\to V'$, $\rho'':W\to V''$ are $\I$-graded $\bfk$-linear isomorphisms. We consider the following diagram
	\[\begin{tikzcd}
		E_{V'}\times E_{V''}&E_1\ar[r,"p_2"]\ar[l,"p_1",swap,]&E_2\ar[r,"p_3"]&E_V
	\end{tikzcd}\]
	where for $(x,W,\rho',\rho'')\in E_1$ and $(x,W)\in E_2$, we set 
	$$p_1(x,W,\rho_1,\rho_2)=(\rho'\cdot x|_{V/W},\rho''\cdot x|_W), \;\; p_2(x,W,\rho_1,\rho_2)=(x,W),\;\;p_3(x,W)=x.$$ 
	
	For a fixed $\I$-graded subspace $W$ in $V$ of dimension $\bd''$ and fixed linear isomorphisms $\rho':V/W\to V'$, $\rho'':W\to V'$, we let $F$ be the subvariety of $E_V$ consisting of $x\in E_V$ such that $W$ is $x$-stable. Then we have the following diagram
	\[\begin{tikzcd}
		E_{V'}\times E_{V''}&F\ar[r,"\iota"]\ar[l,"\kappa",swap]&E_{V}
	\end{tikzcd}\]
	where $\iota$ is the natural inclusion and $\kappa(x)=(\rho'\cdot x|_W,\rho''\cdot x|_{V/W})$ for $x\in F$. Note that $p_1$ is smooth with connected fibers, $p_2$ is a $G_{V'}\times G_{V''}$-principal bundle, and $p_3$ is proper. In \cite{Lus91}, Lusztig defined the induction functor and restriction functor
	\begin{align*}
		\Ind_{V',V''}^V(A\boxtimes B)  &=  (p_3)_! (p_2)_{\flat} (p_1)^* (A\boxtimes B) [d_1-d_2] \\
		\Res_{V',V''}^V(C)             &=  \kappa_! \iota^* (C) [-\langle\bd',\bd''\rangle]
	\end{align*}
	where $A\in \cd(E_{V'})$, $B\in \cd(E_{V''})$, $C\in \cd(E_V)$, $d_1$, $d_2$ are the dimensions of the fibers of $p_1$, $p_2$, respectively. 
	We also note that
	\[
	d_1-d_2 = \sum_{i \in \I} \bd'_i \bd''_i + \sum_{h \in \Omega} \bd'_{s(h)} \bd''_{t(h)}.
	\]
	
	\begin{lemma}[\cite{Lus93}]\label{Lusztig induction on L_nu}
		Let $\nu' \in I_{\bd'}$ and $\nu'' \in I_{\bd''}$. Then the induction functor $\Ind_{V',V''}^V$ satisfies
		\[
		\Ind_{V',V''}^V(L_{\nu'} \boxtimes L_{\nu''}) = L_{\nu'\nu''},
		\]
		where $\nu'\nu''$ denotes the juxtaposition of $\nu'$ and $\nu''$.
	\end{lemma}
	
	\begin{lemma}[\cite{Lus93}]\label{Lusztig restriction on L_nu}
		Let $\nu=(\nu^1,\dots,\nu^m)\in I_\bd$. Then the restriction functor $\Res_{V',V''}^V$ satisfies
		\[
		\Res_{V',V''}^V(L_{\nu}) = \bigoplus (L_{\tau} \boxtimes L_{\omega})[d(\tau,\omega)]
		\]
		where $\tau$ and $\omega$ runs through all $\tau=(\tau^1,\dots,\tau^m) \in I_{\bd'}$, $\omega=(\omega^1,\dots,\omega^m) \in I_{\bd''}$ such that $\nu^l=\tau^l+\omega^l$ for each $l$, and
		\[d(\tau,\omega)=-\sum_{k}\langle\omega^k,\tau^k\rangle_Q-\sum_{k>l}(\langle\tau^k,\omega^l\rangle_Q+\langle\omega^l,\tau^k\rangle_Q).\]
	\end{lemma}
	
	Let $\mathcal{Q}_{V',V''}$ be the subcategory defined as a special case of $\mathcal{Q}_V$, with $Q$ replaced by the disjoint union of two copies (alternatively, it is the subcategory generated by direct sums of tensor products of complexes in $\mathcal{Q}_{V'}$ and $\mathcal{Q}_{V'}$). Then by the previous lemmas, the induction and restriction functors restrict to
	\begin{align*}
		\Ind_{V',V''}^V&: \mathcal{Q}_{V',V''} \to  \mathcal{Q}_{V},\\
		\Res_{V',V''}^V&: \mathcal{Q}_V     \to  \mathcal{Q}_{V',V''}.
	\end{align*}
	
	We define the direct sum $\mathcal{K}_Q=\bigoplus\mathcal{K}_V$, where $\mathcal{K}_V$ is the Grothendieck group of $\mathcal{Q}_V$ for an $\I$-graded vector space $V$ and the direct sum runs over $\I$-graded vector spaces with different dimension vectors. We endow a $\Z[v,v^{-1}]$-module structure on $\mathcal{K}_Q$ by 
	\[v\cdot[L]=[L[-1]].\]
	
	\begin{theorem}[\cite{Lus93}]
		The induction and restriction functors define a bialgebra structure on $\mathcal{K}_Q$ such that $\mathcal{K}_Q$ is a $\Z[v,v^{-1}]$-form of $\bU^+$. Moreover, the class $[L_{n\alpha_i}]$ corresponds to $\mathcal{E}_i^{(n)}$ for $i\in\I$, $n\in\N$.
	\end{theorem}

	In other words, the Grothendieck ring $\mathcal{K}_Q$ can be embedded into $\bU^+$ as a sub-bialgebra and generate $\bU^+$ over $\Q(v^{1/2})$. Note that $\mathcal{K}_Q$ has a $\Z[v,v^{-1}]$-basis $\mathbf{B}=\coprod\mathbf{B}_V$, where $\mathbf{B}_V=\{[L]\mid L\in\mathcal{P}_V\}$. The image of $\mathbf{B}$ under the embedding is called the \textit{canonical basis} of $\bU^+$. We also remark that if $Q$ is a Dynkin quiver, then $\mathbf{B}_V$ consists of $[\IC(\mathcal{O},\overline{\Q}_\ell)]$, where $\mathcal{O}\subset E_V$ is a $G_V$-orbit.

	\begin{example}[\cite{Lus90}]
		\begin{enumerate}
			\item For $\U^+_v(\mathfrak{sl}_2)$, we know that its canonical basis is $\{\ce_1^{(r)}\mid r\in\Z\}$, where $\ce_1^{(r)}=\ce_1^r/[r]!$ is the divided power.
			\item For $\U^+_v(\mathfrak{sl}_3)$, we know that its canonical basis is 
			\begin{align}
				\{\ce_1^{(b)}\ce_2^{(b+c)}\ce_1^{(a)}\mid c\geq a\}\cup \{\ce_2^{(c)}\ce_1^{(a+b)}\ce_2^{(b)}\mid c<a\}.
			\end{align}
		\end{enumerate}    
	\end{example}
	
	\subsection{Canonical basis via Hall algebras}\label{subsec: CB by HA}
	Recall that $\widetilde{\ch}(Q)$ is the generic Hall algebra of $Q$. By using the isomorphism $R^+$ of Theorem~\ref{thm: iso R^pm}, we see $\widetilde{\ch}(Q)$ has an involution $\psi$ given by $v\mapsto v^{-1}$, $\fw_{\alpha_i}\mapsto\fw_{\alpha_i}$. For any $\lambda\in\mathfrak{P}$, we define
	\begin{align}
		E_\lambda=v^{\dim\End(M_q(\lambda))-\dim M_q(\lambda)}\fw_\lambda\in\widetilde{\ch}(Q).
	\end{align}

	Let $\prec$ be the partial order on $\mathfrak{P}$ defined by orbit closure: we say that $\lambda\prec\mu$ if the orbit $\mathfrak{O}_{M_q(\lambda)}$ is contained in the closure of $\mathfrak{O}_{M_q(\mu)}$; see \cite{Lus90} or \cite[\S1.6]{DDPW}. Then following lemma can be deduced using the results in \cite{Lus90}.
	
	\begin{lemma}[\cite{Lus90}]\label{lem:HA bar on E_lambda}
		For any $\lambda\in\mathfrak{P}$, we have
		\begin{align}
			\psi(E_\lambda)\in E_\lambda+\sum_{\mu\prec\lambda}\Z[v,v^{-1}]E_\mu.
		\end{align}
	\end{lemma}
	
	Since for each $\lambda\in\mathfrak{P}$ there are only finitely many $\mu$ such that $\mu\prec\lambda$, we can apply Lusztig's Lemma (\cite[Theorem 1.1]{BZ14}) to extract a bar-invariant basis from $E_\lambda$.
	
	\begin{theorem}
		\label{thm:CB-U+}
		For each $\lambda\in\mathfrak{P}$, there exists a unique element $B_\lambda\in\widetilde{\ch}(Q)$ such that $\psi(B_\lambda)=B_\lambda$ and 
		\[B_\lambda-E_\lambda\in\sum_{\mu\prec\lambda}v^{-1}\Z[v^{-1}]E_\mu.\]
	\end{theorem}
	The basis $\{B_\lambda\mid \lambda\in\mathfrak{P}\}$ is called the canonical basis of $\widetilde{\ch}(Q)$. This elementary construction of canonical basis first appears in \cite{Lus90}. We remark that the upper triangular property in Theorem~\ref{thm:CB-U+} corresponds to the fact that the intersection cohomology complex of $\mathfrak{O}_{M_q(\lambda)}$ is supported on the closure $\overline{\mathfrak{O}_{M_q(\lambda)}}$. In fact, \cite{Lus90} gives the following geometric realization:
	
	\begin{theorem}
		Under the isomorphism $R^+:\U^+\to \widetilde{\ch}(Q)$, the canonical basis of $\U^+$ is sent to the canonical basis of $\widetilde{\ch}(Q)$.
	\end{theorem}
	
	Moreover, using Gabber's purity theorem \cite[5.3.4]{BBD}, \cite{Lus90} proved the following positivity result:
	\[B_\lambda-E_\lambda\in\sum_{\mu\prec\lambda}v^{-1}\N[v^{-1}]E_\mu.\]
	
	Let $Q'$ be a quiver constructed from $Q$ by reversing some arrows.  A Fourier transform $\widetilde{\ch}(\bfk Q)\rightarrow \widetilde{\ch}(\bfk Q')$ is formulated in \cite{Lus90}, and proved in \cite{SV99}. Lusztig \cite[\S13]{Lus90} proved that the Fourier transform preserves the canonical bases of Hall algebras of $Q$ and $Q'$, which convinces the uniqueness of canonical basis of $\U^+$.

	\section{$\imath$Hall algebras and NKS quiver varieties}

	\subsection{$\imath$Quiver algebras}
	\label{subsec:i-quiver}
	
	Let $\varrho$ be an involution of the quiver $Q$. The pair $(Q,\varrho)$ is called an $\imath$quiver. The $\imath$quiver algebra $\Lambda^\imath$ \cite{LW22a} is defined as the quotient $\bfk \ov{Q}/\ov{I}$, where
	\begin{itemize}
		\item[(i)] $\ov{Q}$ is constructed from $Q$ by adding a loop $\varepsilon_i$ at the vertex $i\in Q_0$ if $\btau i=i$, and adding an arrow $\varepsilon_i: i\rightarrow \btau i$ for each $i\in Q_0$ if $\btau i\neq i$;
		\item[(ii)] $\ov{I}$ is generated by
		\begin{itemize}
			\item[(1)] (Nilpotent relations) $\varepsilon_{i}\varepsilon_{\btau i}$ for any $i\in\I$;
			\item[(2)] (Commutative relations) $\varepsilon_i\alpha-\btau(\alpha)\varepsilon_j$ for any arrow $\alpha:j\rightarrow i$ in $Q_1$.
		\end{itemize}
	\end{itemize}

	\begin{example}[$\imath$quiver algebras of diagonal type]
		Let $Q=(Q_0,Q_1)$ be a quiver and $Q^{\dbl} =Q\sqcup  Q^{\diamond}$,  where $Q^{\diamond}$ is an identical copy of $Q$ with a vertex set $\{i^{\diamond} \mid i\in Q_0\}$ and an arrow set $\{ \alpha^{\diamond} \mid \alpha \in Q_1\}$.  We let $\rm{swap}$ be the involution of $Q^{\rm dbl}$ uniquely determined by $\swa(i)=i^\diamond$ for any $i\in Q_0$. We denote the $\imath$quiver algebra $(\Lambda^{\dbl})^\imath$ by $\Lambda$ in the following. 
	\end{example}
	
	By \cite[Corollary 2.12]{LW22a}, $\bfk Q$ is naturally a subalgebra and also a quotient of $\Lambda^\imath$. Viewing $\bfk Q$ as a subalgebra of $\Lambda^{\imath}$, we have a restriction functor
	\[\res: \mod (\Lambda^{\imath})\longrightarrow \mod (\bfk Q).\]
	
	For each $i\in Q_0$, we define a $\bfk$-algebra (which can be viewed as a subalgebra of $\iLa$)
	\begin{align}\label{dfn:Hi}
		\BH _i:=\left\{ \begin{array}{cc}  \bfk[{\textcolor{purple}{\varepsilon_i}}]/({\textcolor{purple}{\varepsilon_i}}^2) & \text{ if }i=\btau i,\\
			\bfk(\xymatrix{i \ar@<0.5ex>@[purple][r]^{\textcolor{purple}{\varepsilon_i}} & \btau i \ar@<0.5ex>@[purple][l]^{\textcolor{purple}{\varepsilon_{\btau i}}}})/( \varepsilon_i\varepsilon_{\btau i},\varepsilon_{\btau i}\varepsilon_i)  &\text{ if } \btau i \neq i .\end{array}\right.
	\end{align}
	Note that $\BH _i=\BH _{\btau i}$ for any $i\in Q_0$.
	
	Recall that $\I_\varrho$ is a chosen set of orbits. We define the following subalgebra of $\Lambda^{\imath}$:
	\begin{equation}  \label{eq:H}
		\BH =\bigoplus_{i\in \ci }\BH _i.
	\end{equation}
	Note that $\BH $ is a radical square zero self-injective algebra. We denote by
	\begin{align}
		\res_\BH :\mod(\iLa)\longrightarrow \mod(\BH )
	\end{align}
	the natural restriction functor. On the other hand, as $\BH $ is a quotient algebra of $\iLa$, every $\BH $-module can be viewed as a $\iLa$-module.
	
	For $i\in Q_0 =\I$, define the indecomposable module over $\BH _i$ (if $i\in \ci$) or over $\BH_{\btau i}$ (if $i\not \in \ci$)
	\begin{align}
		\label{eq:E}
		\E_i =\begin{cases}
			\bfk[\varepsilon_i]/(\varepsilon_i^2), & \text{ if }i=\btau i;
			\\
			\xymatrix{\bfk\ar@<0.5ex>[r]^1 & \bfk\ar@<0.5ex>[l]^0} \text{ on the quiver } \xymatrix{i\ar@<0.5ex>[r]^{\varepsilon_i} & \btau i\ar@<0.5ex>[l]^{\varepsilon_{\btau i}} }, & \text{ if } i\neq \btau i.
		\end{cases}
	\end{align}
	Then $\E_i$, for $i\in Q_0$, can be viewed as a $\iLa$-module and will be called a {\em generalized simple} $\iLa$-module. 
	
	The following important property of the algebra $\Lambda^\imath$ is proved in \cite{LW22a}. It is the starting point of the $\imath$Hall algebra construction explained in the next subsection.
	
	\begin{proposition}[{\cite{LW22a}}]
		The algebra $\Lambda^\imath$ is $1$-Gorenstein.
	\end{proposition}
	
	\subsection{$\imath$Hall algebras}
	\label{subsec:i-Hall}
	
	In \cite{Br13}, Bridgeland used a localization of the Hall algebra of the category of 2-periodic complexes of projective modules in $\rep(\bfk Q)$ to realize $\tU$. This construction was later generalized by Gorsky \cite{Gor18} to define the semi-derived Hall algebras of Frobenius categories. Inspired by their work, \cite{LP21} introduced the concept of semi-derived Ringel-Hall algebras for periodic complexes of hereditary abelian categories. Later this was generalized to $1$-Gorenstein algebras in \cite{LW22a}. In this subsection we briefly recall its construction and basic properties.
	
	We define the twisted Hall algebra $\widetilde{\ch}(\Lambda^\imath)$ to be the $\Q(\sqq^{1/2})$-algebra on the same vector space as $\ch(\Lambda^\imath)$ with twisted multiplication given by
	$$[M]*[N]=\sqq^{\langle \res M,\res N\rangle_{Q}}\sum_{[L]\in\Iso(\mod(\Lambda^\imath))}\frac{|\Ext^1_{\Lambda^\imath}(M,N)_L|}{|\Hom_{\Lambda^\imath}(M,N)|}[L].$$
	
	Let $\cd_{sg}(\Lambda^\imath)
	$ be the singularity category of $\Lambda^\imath$; c.f. \cite{Ha91,Or04}. Set $\mathcal{I}$ to be the subspace of $\widetilde{\ch}(\Lambda^\imath)$ spanned by all differences 
	\begin{align}
		\label{def:I}
		[M]-[N],\quad \text{if $\res_\BH(M)=\res_\BH(N)$ and $M\cong N$ in $\cd_{sg}(\Lambda^\imath)$}.
	\end{align}
	Using the proof of \cite[Proposition 3.8]{LP24} (cf. \cite[Lemmas A11\& A12]{LW22a}), we know that $\mathcal{I}$ is a two-sided ideal of $\widetilde{\ch}(\Lambda^\imath)$, and then we obtain a quotient algebra $\widetilde{\ch}(\Lambda^\imath)/\mathcal{I}$, and denote it by $\widehat{\ch}(\bfk Q,\varrho)$.

	We define a subset $\mathcal{S}$ of $\widehat{\ch}(\bfk Q,\varrho)$ by
	\begin{equation}
		\label{eq:Sca}
		\cs := \{ a[K] \in \widehat{\ch}(\bfk Q,\varrho) \mid a\in \Q(\sqq^{1/2})^\times, \pd K\leq1\}.
	\end{equation}
	By \cite{LW22a,LP24}, the right localization of $\widetilde{\ch}(\Lambda^\imath)/\mathcal{I}$ with respect to $\cs$ exists, which is denoted by $\widetilde{\ch}(\bfk Q,\btau)$ or $\cs\cd\widetilde{\ch}(\Lambda^\imath)$ and called the $\imath$Hall algebra (or the twisted semi-derived Ringel-Hall algebra) of $\mod(\Lambda^\imath)$.
	
	We have $[\K_i]*[\K_j]=[\K_j]*[\K_i]=[\K_i\oplus \K_j]$ in $\widehat{\ch}(\bfk Q,\varrho)$ for any $i,j\in\I$; cf. \cite[Lemma 4.7]{LW22a}. 
	For any $\alpha=(a_i)_{i\in\I}\in\N^{\I}$, we define in $\widehat{\ch}(\bfk Q,\btau)$: 
	$$[\K_\alpha]=[\oplus_{i\in\I}\K_i^{\oplus a_i}]=\prod_{i\in\I}[\K_i]^{a_i}.$$
	Similarly, one can define $[\K_\alpha]$ in $\widetilde{\ch}(\bfk Q,\btau)$ for $\alpha\in\Z^\I$.

	\begin{proposition}[\text{cf. \cite[Proposition 4.9]{LW22a}}]
		\label{basis-iHall}
		
		(1) The algebra $\widetilde{\ch}(\bfk Q,\btau)$ has a (Hall) basis given by
		\begin{align}
			\label{eq:Hallbasis-hat}
			\{[X]*[\K_\alpha]\mid X\in\mod(\bfk Q)\subseteq \mod(\Lambda^\imath), \alpha\in\Z^{\I}\}.
		\end{align}
		
		(2) The algebra $\widehat{\ch}(\bfk Q,\btau)$ has a (Hall) basis given by
		\begin{align}
			\label{eq:Hallbasis-tilde}
			\{[X]*[\K_\alpha]\mid X\in\mod(\bfk Q)\subseteq \mod(\Lambda^\imath), \alpha\in\N^{\I}\}.
		\end{align}
	\end{proposition}
	
	With these basis of $\widetilde{\ch}(\bfk Q,\btau)$ and $\widehat{\ch}(\bfk Q,\btau)$ one can get the Hall algebra realizations of $\hUi$ and $\tUi$. For this one need to verify the relations presented in Proposition~\ref{prop:Serre}, and this is done in \cite[\S 7]{LW22a}. To summarize, one has the following theorem.
	
	\begin{theorem}[\cite{LW22a}]
		\label{lem:Hall-iQG}
		Let $(Q, \btau)$ be a Dynkin $\imath$quiver. Then we have the following isomorphism $\widetilde{\psi}:\tUi|_{v=\sqq}\stackrel{\simeq}{\rightarrow} \widetilde{\ch}(\bfk Q,\btau)$ of $\Q({\sqq^{1/2}})$-algebras, which sends
		\begin{align}
			\label{eq:psi}
			B_i \mapsto \sqq^{-\frac{1}{2}}[S_i], \qquad 
			\tilde{k}_i \mapsto
			\begin{cases}
				[\K_i]&\text{if $\varrho i\neq i$},\\
				\sqq^{-1}[\K_{i}]&\text{if $\varrho i=i$}.
			\end{cases}
		\end{align}    
	\end{theorem}

	\begin{corollary}[Bridgeland's Theorem reformulated] 
		\label{cor:bridgeland}
		Let $Q$ be a Dynkin quiver. Then we have the following isomorphism of $\Q(\sqq^{1/2})$-algebras
		\begin{align*}
			&\widetilde{\Psi}:\tU|_{v=\sqq}\stackrel{\simeq}{\longrightarrow} \widetilde{\ch}(\bfk Q^{\rm dbl},\swa),
			\\
			E_i \mapsto \sqq^{-\frac{1}{2}}[S_i],&\qquad F_i\mapsto \sqq^{-\frac{1}{2}}[S_{i^\diamond}],
			\qquad
			K_i\mapsto [\K_{i^\diamond}],\qquad K_i'\mapsto [\K_i].
		\end{align*}
	\end{corollary}
	
	Similarly, we have the isomorphisms of $\Q(\sqq^{1/2})$-algebras
	\begin{align*}
		\widehat{\Psi}:\hUi|_{v=\sqq}\longrightarrow \widehat{\ch}(\bfk Q,\btau),
		\qquad
		\widehat{\Psi}:\hU|_{v=\sqq}\longrightarrow \widehat{\ch}(\bfk Q^{\rm dbl},\swa).
	\end{align*}
	
	We also need the generic versions of $\widehat{\ch}(\bfk Q,\btau)$ and $\widetilde{\ch}(\bfk Q,\btau)$ for the construction of dual canonical basis. Recall that $\Phi^+$ is the set of positive roots and $\mathfrak{P}:=\mathfrak{P}(Q)$ is the set of functions $\lambda: \Phi^+\rightarrow \N$. For $(\alpha,\nu),(\beta,\mu)\in\Z^\I\times\fp$, there exists a polynomial $\boldsymbol{\varphi}^{\lambda,\gamma}_{\mu,\alpha;\nu,\beta}(v)\in\Z[v,v^{-1}]$ such that
	\[\big([\K_\alpha]\ast[M_q(\mu)]\big)\ast\big([\K_\beta]\ast[M_q(\nu)]\big)=\sum_{\lambda\in\fp,\gamma\in\Z^\I}\boldsymbol{\varphi}^{\lambda,\gamma}_{\mu,\alpha;\nu,\beta}({\sqq})[\K_\gamma]\ast[M_q(\lambda)]
	\]
	in $\widetilde{\ch}(\bfk Q,\btau)$. 
	The generic $\imath$Hall algebra $\tMHg$ 
	is defined to be the  $\Q(v^{1/2})$-space with a basis $\{\K_\alpha*\fu_\lambda\mid \alpha\in\Z^\I,\lambda\in\fp\}$ 
	with multiplication
	\begin{align}
		\label{eq:generic-mult}
		(\K_\alpha*\fu_\mu)*(\K_\beta*\fu_\nu)=\sum_{\lambda\in\fp,\gamma\in\Z^\I}\boldsymbol{\varphi}^{\lambda,\gamma}_{\mu,\alpha;\nu,\beta}(v)\K_\gamma*\fu_\lambda.
	\end{align}
	
	For $\widehat{\ch}(\bfk Q,\btau)$, we can construct its generic version $\widehat{\ch}(Q,\btau)$ with a basis $\{\K_\alpha*\fu_\lambda\mid \alpha\in\N^\I,\lambda\in\fp\}$ similarly.
	
	From Theorem \ref{lem:Hall-iQG}, see also \cite[Theorem 9.8]{LW22a}, we obtain 
	the following isomorphisms of $\Q(v^{1/2})$-algebras
	\begin{alignat*}{2}
		&\widetilde{\Psi}:\tUi\longrightarrow \widetilde{\ch}(Q,\btau),&\qquad &\widehat{\Psi}:\hUi\longrightarrow \widehat{\ch}(Q,\btau),\\
		&\widetilde{\Psi}:\tU\longrightarrow \widetilde{\ch}(Q^{\rm dbl},\swa),&\qquad &\widehat{\Psi}:\hU\longrightarrow \widehat{\ch}(Q^{\rm dbl},\swa).
	\end{alignat*}

	\subsection{NKS quiver varieties}
	\label{subsec:NKS-variety}

	Let $k$ be an algebraically closed field of characteristic zero. 
	Define the \emph{repetition quiver} $\Z Q$ of $Q$ as follows:
	
	$\triangleright$ the set of vertices is $\{(i,p)\in Q_0\times \Z\}$;
	
	$\triangleright$ an arrow $(\alpha,p):(i,p)\rightarrow (j,p)$ and an arrow $(\bar{\alpha},p):(j,p-1)\rightarrow (i,p)$ are given, for any arrow $\alpha:i\rightarrow j$ in $Q$ and any vertex $(i,p)$.
	
	Define the automorphism $\tau$ of $\Z Q$ to be the shift by one unit to the left, i.e., $\tau(i,p)=(i,p-1)$ for all $(i,p)\in Q_0\times \Z$.
	
	By a slight abuse of notation, associated to $\beta:y\rightarrow x$ in $\Z Q$, we denote by $\bar{\beta}$ the arrow that runs from $\tau x\rightarrow y$. 
	Let $k(\Z Q)$ be the \emph{mesh category} of $\Z Q$, that is, the objects are given by the vertices of $\Z Q$ and the morphisms are $k$-linear combinations of paths modulo the ideal spanned by the mesh relations
	$R_x:=\sum\limits_{\alpha:y\rightarrow x}\alpha\bar{\alpha},$
	where the sum runs through all arrows of $\Z Q$ ending at $x$.

	By a theorem of Happel \cite{Ha2}, there is an equivalence
	\begin{align}
		\label{eq:Happelfunctor}
		H: k(\Z Q)\stackrel{\simeq}{\longrightarrow} \Ind \cd_Q,
	\end{align}
	where $\Ind \cd_Q$ denotes the category of indecomposable objects in the bounded derived category $\cd_Q=\cd^b(kQ)$.
	Using this equivalence, we label once and for all the vertices of $\Z Q$ by the isoclasses of indecomposable objects of $\cd_Q$. 
	Note that the action of $\tau$ on $\Z Q$ corresponds to the action of the AR-translation on $\cd_Q$, and this explains the notation $\tau$.

	Let 
	\begin{align}
		\label{eq:C}
		C=
		\{\text{the vertices labeled by }\Sigma^j\tS_i, \text{ for all } i\in Q_0 \text{ and } j\in\Z \}.
	\end{align}
	Let $\Z Q_C$ be the quiver constructed from $\Z Q$ by adding to every vertex $c \in C$ a new object denoted by $\sigma c$ together with arrows $\tau c \rightarrow \sigma c$ and $\sigma c \rightarrow c$; we refer to $\sigma c$, for $c\in C$, as {\em frozen vertices}.

	The \emph{graded NKS category} $\mcr^{\gr}_C$ is defined to be the $k$-linear category with
	\begin{itemize}
		\item[$\triangleright$] objects: the vertices of $\Z Q_C$;
		\item[$\triangleright$] morphisms: $k$-linear combinations of paths in $\Z Q_C$ modulo the ideal spanned by the mesh relations
		$\sum\limits_{\alpha:y\rightarrow x}\alpha\bar{\alpha},$
		where the sum runs through all arrows of $\Z Q_C$ ending at $x\in\Z Q$ (including the new arrow $\sigma x \rightarrow x$ if $x\in C$).
	\end{itemize}
	These categories were formulated in \cite{KS16, Sch19}; here and below NKS stands for Nakajima-Keller-Scherotzke. The work \cite{KS16} was in turn motivated by \cite{Na01,Na04,HL15,LeP13}; also cf. \cite{Qin}.

	\begin{remark}
		\label{ex:cyclic}
		Let $Q$ be a Dynkin quiver. In Nakajima's original construction \cite{Na01,Na04}, the graded Nakajima category $\mcr^{\gr}$ is defined similarly with
		$C$  the set of all vertices of $\Z Q$. We denote by $\cs^{\gr}$ the full subcategory of $\mcr^\gr$ formed by all $\sigma c$, $c\in \Z Q$. 
	\end{remark}

	Let $(Q,\btau)$ be a Dynkin $\imath$quiver.
	Let $F:\cd_Q\rightarrow \cd_Q$ be a triangulated isomorphism induced by an isomorphism $F$ of $\Z Q$; \eqref{eq:Happelfunctor}.
	In this paper, we only consider $F=\Sigma^2$ and $F=F^\imath=\Sigma\circ \widehat{\varrho}$, where $\Sigma$ is the shift functor and $\widehat{\varrho}$ is the derived functor induced by $\varrho$. 
	Note that the subset $C$ in \eqref{eq:C} is $F$-invariant. The isomorphism $F$ of $\Z Q$ can be uniquely lifted to an isomorphism of $\Z Q_C$  by setting $F(\sigma c)=\sigma(Fc)$ for any $c\in C$, and then the functor $F$ of the mesh category $k(\Z Q)$ can be uniquely lifted to $\mcr^{\gr}_C$, which is also denoted by $F$.
	
	Let
	\[
	\mcr=\mcr_{C,F}:=\mcr^{\gr}_C/F,
	\]
	and let $\cs=\cs_{C,F}$ be the full subcategory of $\mcr$ formed by all $\sigma c$ ($c\in C$), following \cite{Sch19}. Then $\mcr$ and $\cs$ are called the \emph{regular NKS category} and the \emph{singular NKS category} of the pair $(F,C)$. The quotient category 
	$$\cp:=\mcr/( \cs ),$$ which is equivalent to $k(\Z Q/F)$ is called the \emph{preprojective NKS category}.  By our assumption, $\cd_Q/F$ is a triangulated category and $\Ind \cd_Q/F\simeq \cp$; cf \cite{Ke05}.

	Let $\cs$ be a singular NKS category, and $\mcr$ its corresponding regular NKS category. An $\mcr$-module $M$ is \emph{stable} (resp. \emph{costable}) if the support of $\soc(M)$ (respectively,  $\Top(M)$) is contained in $\cs_0$. A module is \emph{bistable} if it is both stable and costable.
	
	Let $\bv \in \N^{\mcr_0-\cs_0}$ and $\bw \in \N^{\cs_0}$ be dimension vectors (with finite supports). Denote by $\rep(\bv,\bw,\mcr)$ the variety of $\mcr$-modules of dimension vector $(\bv,\bw)$. Let $\e_x$ denote the characteristic function of $x\in\mcr_0$, which is also viewed as the unit vector supported at $x$. Let $G_\bv:=\prod_{x\in \mcr_0-\cs_0}\GL({\bv(x)}, k)$.
	
	\begin{definition}
		\label{def:NKS}
		The {\em NKS quiver variety}, or simply {\em NKS variety}, $\cm(\bv,\bw)$ is the quotient $\cs t(\bv,\bw)/G_\bv$,  where $\cs t(\bv,\bw)$ is the subset of $\rep(\bv,\bw,\mcr)$ consisting of all stable $\mcr$-modules of dimension vector $(\bv,\bw)$. Define the affine variety
		\begin{equation}
			\label{eq:M0}
			\cm_0(\bv,\bw)=\cm_0(\bv,\bw, \mcr) :=\rep(\bv,\bw,\mcr)\sslash G_\bv
		\end{equation}
		to be the categorical quotient, whose coordinate algebra is $k[\rep(\bv,\bw,\mcr)]^{G_\bv}$.
	\end{definition}
	
	Then $\cm(\bv,\bw)$ is a pure dimensional smooth quasi-projective variety; see \cite[Theorem 3.2]{Sch19}. 
	We define a partial order $\le$ on the set $\N^{\mcr_0}$ as follows:
	\begin{equation}
		\label{eq:leq}
		\text{  $\bv' \leq \bv
			\Leftrightarrow \bv'(x) \le \bv(x), \forall x \in \mcr_0$. Moreover, $\bv' < \bv \Leftrightarrow \bv'\leq \bv$ and $\bv' \neq \bv$.}
	\end{equation}
	Regarding a dimension vector on $\mcr_0-\cs_0$ as a dimension vector on $\mcr_0$ (by extension of zero),
	we obtain by restriction a partial order $\le$ on the set $\N^{\mcr_0-\cs_0}$.
	
	For $\bv', \bv \in \N^{\mcr_0-\cs_0}$ with $\bv'\leq \bv$ and $\bw \in \N^{\cs_0}$,  there is an inclusion
	\[
	\rep(\bv',\bw,\mcr)\longrightarrow \rep(\bv,\bw,\mcr)
	\]
	by taking a direct sum with the semisimple module of dimension vector $\bv-\bv'$. This yields an inclusion
	\[
	\rep(\bv',\bw,\mcr)\sslash G_{\bv'}\longrightarrow \rep(\bv,\bw,\mcr)\sslash G_\bv.
	\]
	Define the affine variety
	\[
	\cm_0(\bw) =\cm_0(\bw, \mcr) :=\colim\limits_{\bv} \cm_0(\bv,\bw)
	\]
	to be the colimit of $\cm_0(\bv,\bw)$ along the inclusions. 
	Then the projection map
	\begin{equation}   \label{eq:pi}
		\pi:\cm(\bv,\bw)\longrightarrow \cm_0(\bv,\bw),
	\end{equation}
	which sends the $G_\bv$-orbit of a stable $\mcr$-module $M$ to the unique closed $G_\bv$-orbit in the closure of $G_\bv M$, is proper; see \cite[Theorem 3.5]{Sch19}. 
	
	Denote by $\cm^{\text{reg}}(\bv,\bw)\subset \cm(\bv,\bw)$ the open subset consisting of the union of closed $G_\bv$-orbits of stable modules, and then
	\[
	\cm_0^{\text{reg}}(\bv,\bw):=\pi(\cm^{\text{reg}}(\bv,\bw))
	\]
	is an open subset of $\cm_0(\bv,\bw)$. \cite[Lemma 3.4]{Sch19} shows that $\cm(\bv,\bw)$ vanishes on all but finitely many dimension vectors $\bv$. Then $\pi$ induces a stratification
	\begin{equation}\label{eq:stratification for M_0}
		\cm_0(\bw)=\bigsqcup_\bv\cm_0^{\text{reg}}(\bv,\bw)
	\end{equation}
	into finitely many smooth locally closed strata $\cm_0^{\text{reg}}(\bv,\bw)$; see \cite[Theorem 3.5]{Sch19} and its proof.
	
	Let $\res:\mod(\mcr)\rightarrow\mod(\cs)$ be the restriction functor. Then $\res$ induces  morphisms of varieties
	\[
	\res:\cm_0(\bw)\longrightarrow \rep(\bw,\cs),\qquad
	\res\circ\pi:\cm(\bv,\bw)\longrightarrow \rep(\bw,\cs).
	\]
	
	Given $\bv \in \N^{\mcr_0-\cs_0}$,  we define a quantum Cartan matrix (cf. \cite{Sch19})
	\begin{align}
		\label{def:Cq}
		\begin{split}
			{\mathcal C}_q \bv:  \mcr_0-\cs_0 & \longrightarrow\Z,
			\\
			({\mathcal C}_q\bv)(x)& =\bv(x)+\bv(\tau x) -\sum_{y\rightarrow x}\bv(y),
			\quad \text{ for }x \in\mcr_0-\cs_0,
		\end{split}
	\end{align}
	where the sum runs over all arrows $y\rightarrow x$ of $\mcr$ with $y\in\mcr_0-\cs_0$.
	Given $\bw \in \N^{\cs_0}$, define a dimension vector
	\[
	\sigma^*\bw:\mcr_0-\cs_0\longrightarrow\N,
	\qquad
	x \mapsto
	\begin{cases}
		\bw(\sigma x), & \text{ if } x\in C,
		\\
		0, & \text{otherwise}.
	\end{cases}
	\]
	Given $\bv \in \N^{\mcr_0-\cs_0}$, define the dimension vector
	\[
	\tau^*\bv: \mcr_0-\cs_0 \longrightarrow \N,
	\qquad
	x \mapsto \bv(\tau x).
	\]
	
	By \cite[Proposition 4.6]{Sch19}, 
	we obtain the following more precise form of a stratification of $\cm_0(\bw)$
	\begin{equation}   \label{eqn:stratification}
		\cm_0(\bw)=\bigsqcup_{\bv:\sigma^*\bw-{\mathcal C}_q\bv\geq0} \cm_0^{\text{reg}}(\bv,\bw).
	\end{equation}
	
	\begin{definition}   \label{def:pair}
		A pair $(\bv,\bw)$ of dimension vectors $\bv \in \N^{\mcr_0-\cs_0}$ and $\bw \in \N^{\cs_0}$ is called $l$\emph{-dominant} if $\sigma^*\bw-{\mathcal C}_q\bv\geq0$, and $(\bv, \bw)$ is called {\em strongly $l$\emph-dominant} if it is $l$\emph-dominant and  $\cm_0^{\reg}(\bv,\bw) \neq \emptyset.$
	\end{definition}
	
	For our application, we will consider the following pairs $(F,C)$:
	
	\begin{definition} [NKS regular/singular categories] 
		\label{def:RS for QG}
		Let $Q$ be a Dynkin quiver. Denote by $\mcr$ and $\cs$ the regular and singular NKS categories associated to the pair $(F=\Sigma^2, C)$, where
		\[C=\{\text{the vertices labeled by }\Sigma^j \tS_i, \text{ for all } i\in Q_0 \text{ and } j\in\Z \}.\]
	\end{definition}

	\begin{definition} [$\imath$NKS regular/singular categories] 
		\label{def:RS for iQG}
		Let $(Q,\btau)$ be a Dynkin $\imath$quiver. Denote by $\mcr^\imath$ and $\cs^\imath$ the regular and singular NKS categories associated to the admissible pair $(F^\imath, C)$, where $F^\imath=\Sigma \widehat{\btau}$ and $C$ is as above.
	\end{definition}

	\subsection{Quantum Grothendieck rings}
	\label{subsec:geometric-iQG}
	
	We assume that $\cs,\mcr$ are the NKS categories considered in the previous subsection. By \cite[Lemma 3.14]{LW21b}, we know that $\cm_0(\bw)\cong \rep(\bw,\cs)$ for any dimension vector $\bw$.
	
	For any two dimension vectors $\alpha,\beta$ of $\cs$, let $V_{\alpha+\beta}$ be a vector space of graded dimension $\alpha+\beta$. Fix a vector subspace $W_0\subset V_{\alpha+\beta}$ of graded dimension $\alpha$, and let
	\[
	F_{\alpha,\beta}:=\{y\in \rep(\alpha+\beta,\cs)\mid y(W_0)\subset W_0\}
	\]
	be the closed subset of $\rep(\alpha+\beta,\cs)$. 
	Then $y \in F_{\alpha,\beta}$ induces a natural linear map $y': V/W_0 \rightarrow V/W_0$, i.e., $y' \in \rep(\beta,\cs)$.
	Hence we obtain the following convolution diagram
	\[
	\rep(\alpha,\cs)\times \rep(\beta,\cs) \stackrel{p}{\longleftarrow} F_{\alpha,\beta}\stackrel{q}{\longrightarrow}\rep(\alpha+\beta,\cs),
	\]
	where $p(y):=(y|_{W_0},y')$ and $q$ is the natural closed embedding.

	Let $\cd_c(\rep(\alpha,\cs))$ be the derived category of constructible sheaves on $\rep(\alpha,\cs)$. We have the following restriction functor (called comultiplication),
	\begin{align*}
		\widetilde{\Delta}^{\alpha+\beta}_{\alpha,\beta}: \cd_c \big(\rep(\alpha+\beta,\cs) \big) \longrightarrow \cd_c \big(\rep(\alpha,\cs) \big)\times \cd_c \big(\rep(\beta,\cs) \big), \quad F\mapsto p_!q^*(F).
	\end{align*}

	Let $\underline{k}_{\cm(\bv,\bw)}$ be the constant sheaf on $\cm(\bv,\bw)$. Denote by $\pi^{\cs}(\bv,\bw)\in \cd_c(\rep(\bw,\cs))$ (or $\pi(\bv,\bw)$ when there is no confusion) the pushforward along $\pi: \cm(\bv,\bw)\rightarrow \cm_0(\bw)\cong \rep(\bw,\cs)$ of $\underline{k}_{\cm(\bv,\bw)}$ with a grading shift:
	\begin{align}
		\label{eq:piPS}
		\pi(\bv,\bw):= \pi_!(\underline{k}_{\cm(\bv,\bw)}) [\dim \cm(\bv,\bw,\mathcal{R})].
	\end{align}
	
	For a strongly $l$-dominant pair $(\bv,\bw)$, let $\cl^{\cs}(\bv,\bw)$ (or $\cl(\bv,\bw)$) be the intersection cohomology (IC for short) complex associated to the stratum $\cm_0^{\reg}(\bv,\bw,\mathcal{R})$ with respect to the trivial local system, that is,
	\begin{align}
		\label{eq:decomp}
		\cl(\bv,\bw)=\IC(\cm_0^{\reg}(\bv,\bw,\mathcal{R})).
	\end{align}
	
	By the transversal slice Theorem \cite{Na01, LW21b}, when $\mathcal{R}$ and $\mathcal{S}$ are the NKS categories described in Definition~\ref{def:RS for QG}, the direct summands appearing in $\pi(\bv,\bw)$ are the shifts of sheaves $\cl(\bv',\bw)$ with $\bv'\leq \bv$; cf. \eqref{eq:leq}, so we have a decomposition
	\begin{equation}
		\label{eqn:decomposition theorem}
		\pi(\bv,\bw)=\sum_{\bv':\sigma^*\bw-{\mathcal C}_q\bv'\geq0,\bv'\leq \bv} a_{\bv,\bv';\bw}(v)\cl(\bv',\bw),
	\end{equation}
	where we denote by $\cf^{\oplus m}[d]$ by $mv^{-d}\cf$ using an indeterminate $v$, for any sheaf $\cf$, $m\in\N$, and $d\in\Z$ (to compare with \cite{LW21b}, set $v=t^{-1}$). Moreover, we have $a_{\bv,\bv';\bw}(v)\in\N[v, v^{-1}]$, $a_{\bv,\bv';\bw}(v^{-1})=a_{\bv,\bv';\bw}(v)$, and $a_{\bv, \bv;\bw}(v)=1$.
	
	We expect the quiver varieties associated to the $\imath$NKS categories $\mathcal{R}^\imath$ and $\mathcal{S}^\imath$ described in Definition~\ref{def:RS for iQG} also satisfy the above property, so we make the following hypothesis (cf. \cite[hypothesis 3.18]{LW21b}):
	
	\begin{hypothesis}\label{iNKS-hypothesis}
		Only IC sheaves with trivial local systems appear in the decomposition of $\pi(\bv,\bw)$.
	\end{hypothesis}
	
	\textbf{Throughout this paper, we shall assume that Hypothesis~\ref{iNKS-hypothesis} holds.}

	For each $\bw\in \N^{\cs_0}$, the Grothendieck group $K_\bw(\mod(\cs))$ is defined as the free abelian group generated by the perverse sheaves $\cl(\bv,\bw)$ appearing in (\ref{eqn:decomposition theorem}), for various $\bv$. It has a distinguished $\Z[v, v^{-1}]$-basis by \eqref{eqn:decomposition theorem}:
	\begin{align}
		\label{eq:bases}
		\{\cl(\bv,\bw) \mid \sigma^*\bw-{\mathcal C}_q\bv\geq0, \cm_0^{\text{reg}}(\bv,\bw)\neq\emptyset\}.
	\end{align}
	Consider the free $\Z[v, v^{-1}]$-module
	\begin{equation}
		\label{eq:Kgr}
		K^{\mathrm{gr}}(\mod(\cs)) := \bigoplus_\bw K_\bw(\mod(\cs)).
	\end{equation}
	Then $\{\widetilde{\Delta}^{\bw}_{\bw_1,\bw_2} \}$ induces a comultiplication $\widetilde{\Delta}$ on $K^{\mathrm{gr}}(\mod(\cs))$.
	
	Introduce a bilinear form $d(\cdot,\cdot)$ on $\N^{\mcr_0-\cs_0}$ by letting
	\begin{equation}\label{definition:d}
		d\big((\bv_1,\bw_1),(\bv_2,\bw_2) \big)=(\sigma^*\bw_1-{\mathcal C}_q\bv_1)\cdot \tau^* \bv_2+\bv_1\cdot \sigma^*\bw_2,
	\end{equation}
	where $\cdot$ denotes the standard inner product, i.e., $\bv' \cdot \bv'' =\sum\limits_{x\in \mcr_0-\cs_0} \bv'(x)\bv''(x)$.
	
	Using the same proof of \cite{VV}, we know the comultiplication $\widetilde{\Delta}$ is coassociative and satisfies
	\begin{align*}
		\widetilde{\Delta}^{\bw}_{\bw_1,\bw_2} \big(\pi(\bv,\bw) \big)=\bigoplus_{\stackrel{\bv_1+\bv_2=\bv}{\bw_1+\bw_2=\bw}} v^{d((\bv_1,\bw_1),(\bv_2,\bw_2))-d((\bv_2,\bw_2),(\bv_1,\bw_1))}\pi(\bv_1,\bw_1)\boxtimes\pi(\bv_2,\bw_2).
	\end{align*}
	
	Denote the graded dual of $K^{\mathrm{gr}}(\mod(\cs))$ by
	\begin{align}
		\label{eq:Kgr2}
		K^{\mathrm{gr}*}(\mod(\cs)) =\bigoplus_{\bw\in \N^{\cs_0}}\Hom_{\Z[v, v^{-1}]} (K_\bw(\mod(\cs)),\Z[v, v^{-1}]).
	\end{align}
	Then as the graded dual of a coalgebra, $K^{\mathrm{gr}*}(\mod(\cs) )$ becomes a $\Z[v, v^{-1}]$-algebra, whose multiplication is denoted by $\ast$.
	Note that $K^{\mathrm{gr}*}(\mod(\cs))$ is a $\N^{\cs_0}$-graded algebra (called the {\em quantum Grothendieck ring}).  It has a distinguished basis
	\begin{align}
		\label{eq:bases dual}
		\{L (\bv,\bw) \mid \sigma^*\bw-{\mathcal C}_q \bv\geq0, \cm_0^{\text{reg}}(\bv,\bw)\neq\emptyset\},
	\end{align}
	dual to the basis in \eqref{eq:bases}.
	
	Under Hypothesis~\ref{iNKS-hypothesis}, each $\pi(\bv,\bw)$ decomposes as shifts of direct sums of $L(\bv,\bw)$, so we can write
	\begin{align}
		\label{equation multiplication}
		L(\bv_1,\bw_1)\ast  L(\bv_2,\bw_2)=\sum_{\bv \geq \bv_1+\bv_2} c_{\bv_1,\bv_2}^\bv(v) L(\bv,\bw_1+\bw_2),
	\end{align}
	where $c_{\bv_1,\bv_2}^\bv(v) \in \N[v, v^{-1}]$.

	\subsubsection{Case of Dynkin $\imath$quivers}
	Let $(Q,\btau)$ be a Dynkin $\imath$quiver and $\cs^\imath$, $\mcr^\imath$ be the $\imath$NKS categories defined in Definition \ref{def:RS for iQG}. By \cite[Lemma 3.14]{LW21b}, we know that $\mod(\cs^\imath)\simeq \mod(\Lambda^\imath)$, and we will identify them in the following. Recall the restriction functor $\res:\mod(\cs^\imath) \rightarrow \mod(kQ)$. It is natural to identify the Grothendieck groups $K_0(\mod(\cs^\imath))$ and $K_0(\mod(kQ))$. Then we define a bilinear form
	$$\langle \bw,\bw'\rangle_{Q,a}=\langle \bw,\bw'\rangle_Q-\langle \bw',\bw\rangle_Q,\quad\forall \bw,\bw'\in\N^{\cs_0^\imath}.$$
	
	Let us fix a square root $v^{1/2}$ of $v$ once for all. As we need a twisting involving $v^{1/2}$ shortly, we shall consider the ring $\cz:=\Z[v^{1/2}, v^{-1/2}]$ and the field $\Q(v^{1/2})$. The coalgebra $K^{\mathrm{gr}}(\mod(\cs^\imath))$ (cf. \eqref{eq:Kgr}) and its graded dual (cf. \eqref{eq:Kgr2}) up to a base change here in the current setting read as follows:
	\begin{gather}
		K^{\mathrm{gr}}(\mod(\cs^\imath)) =\bigoplus_{\bw\in W^{+}} K_\bw(\mod(\cs^\imath)),
		\label{eq:KSi}
		\\
		\hRiZ=\bigoplus_{\bw\in W^{+}} \hR^{\imath}_{\cz,\bw},
		\quad \text{ where } \hR^{\imath}_{\cz,\bw}  =   \Hom_{\cz} \big(K_\bw(\mod(\cs^\imath)),\cz\big).
		\label{eq:tRi}
	\end{gather}
	
	Then $(\hRiZ, \cdot)$ is the $\cz$-algebra corresponding to the coalgebra $K^{\mathrm{gr}}(\mod(\cs^\imath))$ with the {\em twisted} comultiplication
	\begin{equation}
		\label{eq:tw}
		\{\Delta^{\bw}_{\bw^1,\bw^2} :=v^{\frac12\langle \bw^1,\bw^2\rangle_{Q,a}} \widetilde{\Delta}^{\bw}_{\bw^1,\bw^2}\};
	\end{equation}
	in practice, the product sign $\cdot$ is often omitted. We shall need the $\Q(v^{1/2})$-algebra obtained by a base change below:
	\begin{align}
		\label{eq:hRiQZ}
		\hRi  =\Q(v^{1/2}) \otimes_{\cz} \hRiZ.
	\end{align}
	
	For $i\in\I$, we consider dimension vectors $\bv^i,\bw^i$ given by
	\begin{align}
		\label{eq:v^i}
		\bw^{i}&=\e_{\sigma \tS_{ i}}+\e_{\sigma\tS_{\btau i}},
		&&
		\bv^{i}= \sum_{z\in\mcr^\imath_0-\cs^\imath_0} \dim\cp^\imath(\tS_i ,z)\e_z
	\end{align}
	where $\cp^\imath:=\mcr^\imath/(\cs^\imath)$. The algebra $\tRiZ$ is then defined to be the localization of $\tRiZ$ with respect to the multiplicatively closed subset generated by $L(\bv^i,\bw^i)$ ($i\in Q_0$). We also write
	\begin{align}
		\label{eq:tRiQZ}
		\tRi:=\Q(v^{1/2}) \otimes_{\cz} \tRiZ.
	\end{align} 
	
	The following is the main result of \cite{LW21b}.
	
	\begin{theorem}[{\cite[Theorem 6.25]{LW21b}}]
		\label{thm:iQG-sheaf}
		Let $(Q,\varrho)$ be a Dynkin $\imath$quiver. Then there exists an isomorphism of $\Q(v^{\frac{1}{2}})$-algebras $\tilde{\kappa}:\tUi\stackrel{\simeq}{\rightarrow}\tRi$ which  sends
		\begin{align}
			\label{eq:kappa}
			B_i\mapsto L(0,\e_{\sigma\tS_i}),  \qquad
			\tilde{k}_i \mapsto \begin{cases}
				v^{-1}L(\bv^{\varrho i},\bw^i)&\textit{if $\varrho i=i$},\\
				L(\bv^{i},\bw^i)&\textit{if $\varrho i\neq i$}.
			\end{cases}
		\end{align}
	\end{theorem}

	\subsubsection{Case of double framed quivers}
	Let $\cs,\mcr$ be the NKS categories defined in Definition \ref{def:RS for QG}. Using the restriction functor $\res:\mod(\cs) \rightarrow \mod(kQ^{\rm dbl})$, we may identify the Grothendieck groups $K_0(\mod(\cs))$ and $K_0(\mod(kQ^{\rm dbl}))$. Let $\langle\cdot,\cdot\rangle_{Q^{\dbl}}$ be the Euler form of $kQ^{\dbl}$. We define the bilinear form $\langle\cdot,\cdot\rangle_{Q^{\rm dbl},a}$ as follows: for any dimension vectors $\bw,\bw' \in \N^{\cs_0}$, let
	\begin{align}
		\langle \bw,\bw'\rangle_{Q^{\rm dbl},a}&= \langle \bw,\bw'\rangle_{Q^{\rm dbl}}-\langle \bw',\bw\rangle_{Q^{\rm dbl}}.\label{eqn: antisymmetric bilinear form}
	\end{align}
	
	The coalgebra $K^{\mathrm{gr}}(\mod(\cs))$ (cf. \eqref{eq:Kgr}) and its graded dual (cf. \eqref{eq:Kgr2}) up to a base change here in the current setting read as follows:
	\begin{gather}
		K^{\mathrm{gr}}(\mod(\cs)) =\bigoplus_{\bw\in W^{+}+W^{-}} K_\bw(\mod(\cs)),
		\label{eq:KS}
		\\
		\hRZ=\bigoplus_{\bw\in W^{+}+W^{-}} \hR_{\cz,\bw},
		\quad \text{where } \hR_{\cz,\bw}  =   \Hom_{\cz} \big(K_\bw(\mod(\cs)),\cz\big).
		\label{eq:tR}
	\end{gather}
	
	Then $(\hRZ, \cdot)$ is the $\Z[v^{1/2}, v^{-1/2}]$-algebra corresponding to the coalgebra $K^{\mathrm{gr}}(\mod(\cs))$ with the {\em twisted} comultiplication
	\begin{equation}
		\label{eq:tw}
		\Delta^{\bw}_{\bw^1,\bw^2} :=v^{\frac12\langle \bw^1,\bw^2\rangle_{Q^{\dbl},a}} \widetilde{\Delta}^{\bw}_{\bw^1,\bw^2}.
	\end{equation}
	We also need the $\Q(v^{1/2})$-algebra obtained by a base change:
	\begin{align}
		\label{eq:hRQZ}
		\hR  =\Q(v^{1/2}) \otimes_{\cz} \hRZ.
	\end{align}
	
	Again for $i\in\I$, we consider dimension vector $\bv^i$, $\bv^{\Sigma i}$ and $\bw^i$ given by
	\begin{align}
		\bw^{i}&=\e_{\sigma \tS_{ i}}+\e_{\sigma\Sigma\tS_{i}},
		\qquad
		\bv^{i}= \sum_{z\in\mcr_0-\cs_0} \dim\cp(\tS_i ,z)\e_z,\qquad
		\bv^{\Sigma i}=\Sigma^*\bv^i,
	\end{align}
	and set
	\begin{align}
		&W^{+}=\bigoplus_{x\in\{\tS_i,i\in Q_0\}}\N \e_{\sigma x},&&
		V^{+}=\bigoplus_{x\in\Ind \mod(kQ),\, x\text{ is not injective}} \N \e_x,\\
		&W^{-}=\Sigma^*W^{+},&&
		V^{-}=\Sigma^*V^{+},\\
		&W^{0}=\bigoplus_{i\in Q_0} \N \bw^{i},&&
		V^{0}=\bigoplus_{i\in Q_0} \N \bv^{i}\oplus \N \bv^{\Sigma i}.
	\end{align}
	Let $\hR^+$ be the submodule of $\hR$ generated by $L(\bv,\bw)$, $(\bv,\bw)\in (V^+,W^+)$. It also becomes an algebra under $\Delta^\bw_{\bw_1+\bw_2}$, and the inclusion $\hR^+\hookrightarrow\hR$ is an algebra homomorphism. Similarly we can define the subalgebra $\hR^-$.
	
	Finally, we let $\tRZ$ be the localization of $\tRZ$ with respect to the multiplicatively closed subset generated by $L(\bv^i,\bw^i)$ and $L(\bv^{\Sigma i},\bw^i)$ ($i\in Q_0$), and
	\begin{align}
		\label{eq:tRQZ}
		\tR: =\Q(v^{1/2}) \otimes_{\cz} \tRZ.
	\end{align} 
	
	The following theorem recovers the main result of \cite{Qin}.
	
	\begin{theorem}
		Let $(Q,\varrho)$ be a Dynkin $\imath$quiver. Then there exists an isomorphism of $\Q(v^{\frac{1}{2}})$-algebras $\tilde{\kappa}:\tU\stackrel{\simeq}{\rightarrow}\tR$ which sends
		\begin{align*}
			E_i\mapsto L(0,\e_{\sigma\tS_i}),\quad F_i\mapsto L(0,\e_{\sigma \Sigma \ts_i}),\quad K_i\mapsto L(\bv^i,\bw^i),\quad K_i'\mapsto L(\bv^{\Sigma i},\bw^i).
		\end{align*}
	\end{theorem}
	
	Note that by restricting the isomorphism $\tilde{\kappa}:\tU\stackrel{\simeq}{\rightarrow}\tR$ to $\U^\pm$, we get $\Q(v^{1/2})$-algebra isomorphisms
	\begin{align}
		&\kappa^\pm:\U^\pm\longrightarrow \tR^{\pm}, 
		\\
		E_i\mapsto &L(0,\e_{\sigma\tS_i}),\quad F_i\mapsto L(0,\e_{\sigma \Sigma \ts_i}).
	\end{align}
	
	Recall the generators $\ce_i$, $\cf_i$ in \eqref{eq:Udj-gen}. Let $\mathbf{B}^-$ be Lusztig's canonical basis of $\U^-$ (via the generators $\cf_i$ ($i\in\I$)). Following \cite{Ka91}, we define a Hopf pairing on $\U^-\otimes\U^+$ by
	\begin{align} 
		\label{eq:hopf}(F_i,E_j)_{K}=\delta_{ij}(v-v^{-1}).
	\end{align}
	
	For $b\in\mathbf{B}^-$, we denote by $\delta_b\in\U^+$ the dual basis of $b$ under the pairing $(\cdot,\cdot)_K$. We define a norm function $N:\Z^\I\rightarrow\Z$ by
	\[N(\alpha)=\frac{1}{2}(\alpha,\alpha)_Q-\eta(\alpha).\]
	where $\eta:\Z^\I\rightarrow\Z$ is the augmentation map defined by $\eta(\sum_ia_i\alpha_i)=\sum_ia_i$ (cf. \cite{HL15}). The rescaled dual canonical basis of $\U^+$ is then defined to be
	\[\widetilde{\mathbf{B}}^+:=\{v^{\frac{1}{2}N(-\deg(b))}\delta_b\mid b\in\mathbf{B}^-\}.\]
	The rescaled dual canonical basis of $\U^-$ is defined as $\varphi(\widetilde{\mathbf{B}}^+)$, where $\varphi$ is the isomorphism $\U^+\rightarrow\U^-$ defined by $E_i\mapsto F_i$.

	\begin{example}
		For $\tU=\tU_v(\mathfrak{sl}_2)$, we have $\widetilde{\mathbf{B}}^+=\{E_1^n\mid n\in\N\}$ and $\widetilde{\mathbf{B}}^-=\{F_1^n\mid n\in\N\}$.
	\end{example}

	\begin{theorem}[{\cite{HL15, Qin}}]\label{dCB of U^+ by L}
		Under the isomorphism $\kappa^+:\U^+\rightarrow\tR^+$, the basis $\{L(\bv,\bw)\mid (\bv,\bw)\in(V^+,W^+)\}$ gets identified with the rescaled dual canonical basis of $\U^+$.
	\end{theorem}
	
	\section{Dual canonical bases}

	\subsection{Dual canonical bases of Hall algebras}\label{dCB of HA subsec}
	
	For a Dynkin quiver $Q$, we let $\widetilde{\ch}(Q)$ be the generic Hall algebra of $Q$. From \eqref{eq:mult}, we know that it also has a basis $\{\fu_\lambda\mid\lambda\in\mathfrak{P}\}$, where $\fu_\lambda=a_\lambda\fw_\lambda$ and $a_\lambda$ is the polynomial such that $a_\lambda(\sqq)=a_{M_q(\lambda)}$. The multiplication formula then becomes
	\[\fu_\mu\cdot\fu_\nu=\sum_{\lambda}g^{\lambda}_{\mu,\nu}(v)\fu_\lambda\]
	where $g^{\lambda}_{\mu,\nu}(v)\in\Z[v,v^{-1}]$ is the polynomial such that
	\[
	[M_q(\mu)]\cdot[M_q(\lambda)]=\sum_{\lambda\in\mathfrak{P}}g^{\lambda}_{\mu,\nu}(\sqq)[M_q(\lambda)].
	\]
	In this subsection we will extend the coefficient ring to $\Q(v^{\frac{1}{2}})$. The following result follows directly from Theorem~\ref{thm: iso R^pm} (note the coefficient on $\fu_{\alpha_i}$):
	
	\begin{proposition}[\cite{Rin90}]
		There exists an isomorphism of $\Q(v^{1/2})$-algebras
		\begin{align*}
			\Psi^+:\U^+\longrightarrow \widetilde{\ch}(Q),\quad E_i\mapsto v^{-\frac{1}{2}}\fu_{\alpha_i},\quad\forall i\in\I.
		\end{align*}
	\end{proposition}
	
	Similarly, there exists an isomorphism of $\Q(v^{1/2})$-algebras 
	\[
	\Psi^-:\U^-\longrightarrow \widetilde{\ch}(Q),\quad F_i\mapsto v^{-\frac{1}{2}}\fu_{\alpha_i},\quad\forall i\in\I.
	\] 
	The Hopf pairing in \eqref{eq:hopf} can then be transferred to $\widetilde{\ch}(\bfk Q)\otimes\widetilde{\ch}(\bfk Q)$, which is given by
	\begin{align}
		\label{eq:hopf-Hall}
		(M,N)_K=\delta_{M,N}|\aut(M)|.
	\end{align}
	We can also define a similar paring on $\widetilde{\ch}(Q)\otimes \widetilde{\ch}(Q)$.
	
	Note that there is an algebra isomorphism 
	\[\Omega^+:\widetilde{\ch}(Q)\longrightarrow \tR^+,\quad \fu_{\alpha_i}\mapsto v^{\frac{1}{2}}L(0,\e_{\sigma\tS_i}),\quad\forall i\in\I,\]
	such that $\Omega^+\circ\kappa^+=\Psi^+$.
	The inverse image of $\{L(\bv,\bw)\mid(\bv,\bw)\in(V^+,W^+)\}$ under $\Omega^+$ will be called the dual canonical basis of $\widetilde{\ch}(Q)$. In view of Theorem~\ref{dCB of U^+ by L}, it is also the inverse image of $\widetilde{\mathbf{B}}^+$ under $\Psi^+$.
	
	\begin{example}
		For $\tU_v(\mathfrak{sl}_2)$, the dual canonical basis of $\widetilde{\ch}(Q)$ is $\{v^{-\frac{m^2}{2}}\fu_{m\alpha_1}\mid m\in\N\}$.
	\end{example}
	
	Recall from Proposition~\ref{QG bar-involution def} that there is a bar-involution on $\U^+$ defined by $v^{\frac{1}{2}}\mapsto v^{-\frac{1}{2}}$ and $E_i\mapsto E_i$, which is an anti-automorphism. Using the isomorphism $\psi^+$, this can be transferred to the generic Hall algebra $\widetilde{\ch}(Q)$: it is characterized by $\ov{\fu_{\alpha_i}}=v^{-1}\fu_{\alpha_i}$. 
	
	For $\lambda\in\mathfrak{P}$, we define a rescaled Hall basis of $\widetilde{\ch}(Q)$ by
	\begin{equation}\label{eq:HA element U_lambda}
		\mathfrak{U}_\lambda=v^{-\dim\End_{\bfk Q}(M_q(\lambda))+\frac{1}{2}\langle M_q(\lambda),M_q(\lambda)\rangle_Q}\fu_\lambda.
	\end{equation}
	
	Then the following analogue of Lemma~\ref{lem:HA bar on E_lambda} is straightforward to verify:
	\begin{lemma}
		\label{lem: HA bar of U_lambda}
		For each $\lambda\in\mathfrak{P}$, one can write
		\[\ov{\mathfrak{U}_\lambda}-\mathfrak{U}_\lambda\in\sum_{\lambda\prec\mu}\Z[v,v^{-1}]\mathfrak{U}_\mu.\]
	\end{lemma}
	
	Using the same idea as for the canonical basis of $\widetilde{\ch}(Q)$, we can apply Lusztig's Lemma (\cite[Theorem 1.1]{BZ14}) to produce a bar-invariant basis from $\mathfrak{U}_\lambda$. As we shall see, this construction actually gives the dual canonical basis of $\widetilde{\ch}(Q)$.
	
	\begin{theorem}
		\label{thm:dualCB-U+}
		For each $\lambda\in\mathfrak{P}$, there exists a unique element $\mathfrak{C}_\lambda\in\widetilde{\ch}(Q)$ such that $\ov{\mathfrak{C}_\lambda}=\mathfrak{C}_\lambda$ and 
		\[\mathfrak{C}_\lambda-\mathfrak{U}_\lambda\in\sum_{\mu\in\mathfrak{P}}v^{-1}\Z[v^{-1}]\mathfrak{U}_\mu.\]
		Moreover, $\mathfrak{C}_\lambda$ satisfies
		\[\mathfrak{C}_\lambda-\mathfrak{U}_\lambda\in\sum_{\lambda\prec\mu}v^{-1}\Z[v^{-1}]\mathfrak{U}_\mu.\]
	\end{theorem}
	
	We identify the basis $\{\mathfrak{C}_\lambda\mid\lambda\in\fp\}$ of $\widetilde{\ch}(Q)$, let us first recall the construction of canonical basis given in \ref{subsec: CB by HA}. For each $\lambda\in\fp$, we have defined a basis element of $\widetilde{\ch}(Q)$ by
	\[
	E_\lambda=v^{\dim \End(M_q(\lambda))-\frac{1}{2}\dim M_q(\lambda)}\fw_\lambda,
	\]
	The corresponding canonical basis $B_\lambda\in \U^+$ then satisfies
	\[
	B_\lambda-E_\lambda\in\sum_{\mu\prec\lambda}v^{-1}\Z[v^{-1}]E_\mu.
	\]
	
	Let $E_\lambda^*$ (resp. $B_\lambda^*$) be the dual basis of $E_\lambda$ (resp. $B_\lambda$) with respect to the paring \eqref{eq:hopf-Hall}. We also introduce a function $N:\Z^\I\to\Z$ defined by
	\[N(\alpha)=\frac{1}{2}(\alpha,\alpha)_Q-\eta(\alpha).\]
	where $\eta:\Z^\I\rightarrow\Z$ is the augmentation map defined by $\eta(\sum_ia_i\alpha_i)=\sum_ia_i$ (cf. \cite{HL15}). Then one can check that
	\begin{align}
		\label{eq:TE}
		\mathfrak{U}_\lambda=v^{\frac{1}{2}n_\lambda}E_\lambda^*,
	\end{align}
	where $n_\lambda=N(\sum_k\lambda(\beta_k)\beta_k)$. Moreover,
	\begin{align*}
		B_\lambda^*\in E_\lambda^*+\sum_{\lambda\prec \mu}v^{-1}\Z[v^{-1}]E_\mu^*,\quad \ov{B_\lambda^*}=v^{n_\lambda}B_\lambda^*.
	\end{align*}
	In view of \eqref{eq:TE} and Theorem~\ref{thm:dualCB-U+}, this means $v^{\frac{1}{2}n_\lambda}B_\lambda^*=\mathfrak{C}_\lambda$ (note that $\lambda\prec\mu$ implies $\dimv M_q(\lambda)=\dimv M_q(\mu)$, which means $n_\lambda=n_\mu$). The following proposition then follows Theorem~\ref{dCB of U^+ by L}.
	
	\begin{proposition}\label{dCB of HA is L}
		We have $\Omega^+(\mathfrak{C}_\lambda)=L(\bv_\lambda,\bw_\lambda)$ for any $\lambda\in\fp$, where $(\bv_\lambda,\bw_\lambda)\in(V^+,W^+)$ is such that $\sigma^*\bw_\lambda-\mathcal{C}_q\bv_\lambda=\lambda$.
	\end{proposition}
	
	It is the construction of the dual canonical basis $\{\mathfrak{C}_\lambda\mid\lambda\in\fp\}$ that will be generalize to $\imath$Hall algebras.
	
	\subsection{Dual canonical basis of $\imath$Hall algebras}
	\label{subsec:dCB-iHall}
	
	Inspired by the results of \S\ref{dCB of HA subsec}, we will also construct a bar-invariant basis of $\widehat{\ch}(Q,\varrho)$. Recall that there is a natural inclusion of $\Q(v^{\frac{1}{2}})$-vector spaces
	\begin{equation}\label{eq:inclusion of HA to iHA}
		\iota:\widetilde{\ch}(Q)\hookrightarrow \widehat{\ch}(Q,\varrho),\quad \fu_\lambda\mapsto\fu_\lambda,\forall\lambda\in\mathfrak{P}.
	\end{equation}
	The element $\mathfrak{U}_\lambda$ defined by \eqref{eq:HA element U_lambda} can therefore be regarded as an element in $\widehat{\ch}(Q,\varrho)$. Let $\widehat{\ct}(Q,\varrho)$ be the subalgebra generated by $\K_\alpha$, $\alpha\in\N^\I$. Then $\widehat{\ch}(Q,\varrho)$ has a basis given by 
	\[
	\{\K_\alpha\ast\fu_\lambda\mid\alpha\in\N^\I,\lambda\in\mathfrak{P}\}.
	\]
	
	The bar-involution of $\widehat{\ch}(Q,\varrho)$ can be defined similarly to $\widehat{\ch}(Q,\varrho)$ by setting
	\[\ov{\fu_{\alpha_i}}=v^{-1}\fu_{\alpha_i},\quad \ov{\K_{\alpha_i}}=\K_{\alpha_i}.\]
	Following \cite{BG17}, we define an action of $\widehat{\ct}(Q,\varrho)$ on $\widehat{\ch}(Q,\varrho)$ by
	\[\K_\alpha\diamond\fu_\lambda=v^{\frac{1}{2}(\alpha-\varrho\alpha,\,\dimv{M_q(\lambda)})_Q}\K_\alpha\ast\fu_\lambda.\]
	This definition has the advantage that for any $\alpha\in\N^\I$ and $\lambda\in\mathfrak{P}$,
	\begin{equation}\label{eq:diamond action and bar}
		\ov{\K_\alpha\diamond \fu_\lambda}=\K_\alpha\diamond\ov{\fu_\lambda}.
	\end{equation}
	
	Similar to the case of Hall algebras, we define a partial order on $\N^\I\times\mathfrak{P}$: we say $(\alpha,\lambda)\prec(\beta,\mu)$ if $\alpha+\btau(\alpha)+\dimv M_q(\lambda)=\beta+\btau(\beta)+\dimv M_q(\mu)$ and either $\alpha\prec\beta$ (i.e. $\alpha\neq\beta$ and $\beta-\alpha\in\N^\I$) or $\alpha=\beta$ and $\lambda\prec\mu$. The following result is an analogue of Lemma~\ref{lem: HA bar of U_lambda} for $\imath$Hall algebras.
	
	\begin{lemma}[{\cite{LP25}}]\label{iHA bar of H_lambda}
		For $\alpha\in\N^\I$ and $\lambda\in\mathfrak{P}$, we have 
		\[\ov{\K_\alpha\diamond \mathfrak{U}_\lambda}-\K_\alpha\diamond \mathfrak{U}_\lambda\in\sum_{(\alpha,\lambda)\prec(\beta,\mu)}\Z[v,v^{-1}]\cdot \K_\beta\diamond \mathfrak{U}_\mu.\]
	\end{lemma}
	
	Since for a fixed pair $(\alpha,\lambda)$ there are only finitely many pairs $(\beta,\mu)$ such that $(\alpha,\lambda)\prec(\beta,\mu)$, Lusztig's Lemma \cite[Theorem 1.1]{BZ14} is applicable and we can construct a bar-invariant basis for $\widehat{\ch}(Q,\varrho)$ from $\{\K_\alpha\diamond \mathfrak{U}_\lambda\mid \alpha\in\N^\I,\lambda\in\fp\}$. We also note the $\diamond$-action preserves this basis, more precisely, the following theorem is valid. 
	
	\begin{theorem}[{\cite{LP25}}]\label{iHA dCB theorem}
		For each $\alpha\in\N^\I$ and $\lambda\in\mathfrak{P}$, there exists a unique element $\mathfrak{L}_{\alpha,\lambda}\in\widehat{\ch}(Q,\varrho)$ such that $\ov{\mathfrak{L}_{\alpha,\lambda}}=\mathfrak{L}_{\alpha,\lambda}$ and
		\[
		\mathfrak{L}_{\alpha,\lambda}-\K_\alpha\diamond \mathfrak{U}_\lambda\in\sum_{(\beta,\mu)}v^{-1}\Z[v^{-1}]\cdot \K_\beta\diamond \mathfrak{U}_\mu.
		\]
		Moreover, $\mathfrak{L}_{\alpha,\lambda}$ satisfies 
		\[
		\mathfrak{L}_{\alpha,\lambda}-\K_\alpha\diamond \mathfrak{U}_\lambda\in\sum_{(\alpha,\lambda)\prec(\beta,\mu)}v^{-1}\Z[v^{-1}]\cdot \K_\beta\diamond \mathfrak{U}_\mu,
		\]
		and $\mathfrak{L}_{\alpha,\lambda}=\K_\alpha\diamond \mathfrak{L}_{0,\lambda}$.
	\end{theorem}

	Because of the last property, we often write $\mathfrak{L}_\lambda:=\mathfrak{L}_{0,\lambda}$ for $\lambda\in\mathfrak{P}$ and use $\{\K_\alpha\diamond \mathfrak{L}_\lambda\mid \alpha\in\N^\I,\lambda\in\mathfrak{P}\}$ to denote the basis constructed in Theorem~\ref{iHA dCB theorem}. This is called the dual canonical basis of $\widehat{\ch}(Q,\varrho)$. The dual canonical basis of $\widetilde{\ch}(Q,\varrho)$ is then defined to be $\{\K_\alpha\diamond \mathfrak{L}_\lambda\mid \alpha\in\Z^\I,\lambda\in\mathfrak{P}\}$.

	Using the isomorphisms given in Theorem \ref{thm:iQG-sheaf} and Lemma \ref{lem:Hall-iQG}, we obtain two isomorphisms of $\Q(v^{1/2})$-algebras:
	\begin{align}
		\widetilde{\Omega}:&\widetilde{\ch}(Q,\btau)\longrightarrow \tR^{\imath},\qquad \widehat{\Omega}:\widehat{\ch}(Q,\btau)\longrightarrow \hRi,
	\end{align}
	defined by $
	\fu_{\alpha_i}\mapsto v^{\frac{1}{2}}L(0,1_{\sigma S_i})$ and $\K_{\alpha_i}\mapsto L(\mathbf{v}^{\varrho i},\mathbf{w}^i)$.

	\begin{proposition}[{\cite[Proposition 6.5]{LP26a}}]
		\label{iHA dCB compare}
		For $\alpha\in\Z^\I$ and $\lambda\in\mathfrak{P}$, the dual canonical basis $\K_\alpha\diamond \mathfrak{L}_\lambda$ of $\widetilde{\mathcal{H}}(Q,\varrho)$ is mapped to $L(\bv_\alpha+\bv_\lambda,\bw_\alpha+\bw_\lambda)$ under the isomorphism $\widetilde{\Omega}:\widetilde{\ch}(Q,\varrho)\rightarrow\tRi$, where $(\bv_\lambda,\bw_\lambda)$ is such that $\sigma^*\bw_\lambda-\cc_q\bv_\lambda=\lambda$ and $\widetilde{\Omega}(\K_\alpha)=L(\bv_\alpha,\bw_\alpha)$.
	\end{proposition}
	
	\begin{remark}
		In order to prove Proposition~\ref{iHA dCB compare}, the authors of \cite{LP26a} had to introduce the open subvariety of $\rep(\bw,\mathcal{S})$ consisting of Gorenstein projective modules and consider the corresponding quantum Grothendieck ring. In the context of Hall algebras, this means we pass from the semi-derived Hall algebra $\widetilde{\ch}(Q,\varrho)=\cs\cd\widetilde{\ch}(\Lambda^\imath)$ to $\mathcal{H}(\Gproj(\Lambda^\imath))$, the Hall algebra of the category of Gorenstein projective $\Lambda^\imath$-modules (which is a subalgebra of $\cs\cd\widetilde{\ch}(\Lambda^\imath)$, and they are isomorphic after doing  localization with respect to projective $\Lambda^\imath$-modules in $\mathcal{H}(\Gproj(\Lambda^\imath))$; see \cite{Gor18,LW22a}). The construction of dual canonical basis can also be stated for $\mathcal{H}(\Gproj(\Lambda^\imath))$. Then using a geometric approach similar to \cite{Lus90}, one can prove the analogue of Proposition~\ref{iHA dCB compare} for $\mathcal{H}(\Gproj(\Lambda^\imath))$, which will in turn imply the corresponding result for $\widetilde{\ch}(Q,\varrho)$.
	\end{remark}
	
	From the result of Proposition~\ref{iHA dCB compare} one can deduce various positivity results for the dual canonical basis of $\widetilde{\ch}(Q,\btau)$. For example, we have the following:
	
	\begin{corollary}[{\cite[Theorem 6.6]{LP25}}]
		For $\alpha\in\N^\I$ and $\lambda\in\fp$, in $\widetilde{\ch}(Q,\varrho)$ we have 
		\[\K_\alpha\diamond\mathfrak{U}_\lambda\in \K_\alpha\diamond\mathfrak{L}_\lambda+\sum_{(\alpha,\lambda)\prec(\beta,\mu)}v^{-1}\N[v^{-1}]\K_\beta\diamond \mathfrak{L}_\mu.\] 
	\end{corollary}

	\subsection{Reflection functors}\label{subsec:iHA reflection}
	Let $\ell\in Q_0$ be a sink. Denote by $s_\ell Q$ the quiver obtained from $Q$ by reversing all arrows ending at $\ell$. Let $\mathscr{R}_\ell^+:\mod(\bfk Q)\rightarrow \mod(\bfk s_\ell Q)$ be the BGP reflection. 
	
	For an $\imath$quiver $(Q,\btau)$, define for any $\ell\in Q_0$
	\begin{align*}
		r_\ell=\begin{cases}
			s_\ell s_{\btau \ell} & \text{ if }\btau \ell\neq \ell,
			\\
			s_\ell &\text{ if }\btau \ell=\ell.
		\end{cases}
	\end{align*}
	In this way, we define $r_\ell Q$ for any sink $\ell\in Q_0$, and then $(r_\ell Q,\varrho)$ is also an $\imath$quiver. Denote by $r_\ell\Lambda^\imath$ its $\imath$quiver algebra.
	In \cite{LW21a}, we define a BGP type reflection 
	$F_\ell^+:\mod(\Lambda^\imath)\rightarrow\mod(r_\ell \Lambda^\imath)$, which restricts to $\mod(\bfk Q)$ is 
	\begin{align*}
		\begin{cases}
			\mathscr{R}_\ell^+  :\mod(\bfk Q)\rightarrow\mod(\bfk r_\ell Q)&\text{ if }\btau\ell=\ell,
			\\
			\mathscr{R}_\ell^+ \mathscr{R}_{\varrho\ell}^+=\mathscr{R}_{\varrho\ell}^+ \mathscr{R}_\ell^+ :\mod(\bfk Q)\rightarrow\mod(\bfk r_\ell Q)&\text{ if }\btau\ell\neq\ell.
		\end{cases} 
	\end{align*}

	For any $i\in \I$, denote by $S_i$ (respectively, $S_i'$) the simple $\bfk Q$-module (respectively, $\bfk r_\ell Q$-module), denote by $\K_i$ (respectively, $\K_i'$) the
	generalized simple $\iLa$-module (respectively, $\bs_\ell\iLa$-module). We similarly define $[\K_\alpha],[\K_\beta']$ for $\alpha\in K_0(\mod(\bfk Q))$, $\beta\in K_0(\mod(\bfk(r_\ell Q)))$
	in the $\imath$Hall algebras (where $\ell$ is a sink of $Q$).
	
	Recall the root lattice $\Z^{\I}=\Z\alpha_1\oplus\cdots\oplus\Z\alpha_n$, and we have an isomorphism of abelian groups $\Z^\I\rightarrow K_0(\mod(\bfk Q))$, $\alpha_i\mapsto \widehat{S_i}$. This isomorphism induces the action of the reflection $s_i$ on $K_0(\mod(\bfk Q))$.
	Thus for $\alpha\in K_0(\mod (\bfk Q))$ and $i\in\I$, $[\K_{s_i\alpha}] \in \iH(\bfk Q,\btau)$ is well defined. Similarly, we have $[\K'_{s_i\alpha}]\in \iH(\bfk r_\ell Q,\btau)$.
	
	By \cite[Proposition 4.4]{LW21a}, there exists an isomorphism 
	\[\Gamma_\ell:\iH(\bfk Q,\btau)\rightarrow\iH(\bfk r_\ell Q,\btau),\]
	which induces an isomorphism
	\[
	\Gamma_\ell: \iH( Q,\btau)\rightarrow\iH( r_\ell Q,\btau).
	\]
	Moreover, we have the following commutative diagram:
	\begin{equation}\label{eq:defT}
		\begin{tikzcd}
			\tUi  \ar[r,"\tTT_\ell"] \ar[d,swap,"\widetilde{\Psi}_Q"] & \tUi \ar[d,"\widetilde{\Psi}_{r_\ell Q}"]\\
			\tMHg \ar[r,"\Gamma_\ell"] &  \widetilde{\ch}(r_\ell Q,\btau)
		\end{tikzcd}
	\end{equation}
	
	Since $\widetilde{\TT}_\ell$ commutes with the bar-involution of $\tUi$ by Lemma \ref{lem:QGbraid-bar}, we also have $\Gamma_\ell$ commutes with the bar-involution of $\widetilde{\ch}(Q,\varrho)$.

	For $M\in\mod(\bfk Q)$, let us write
	\[\mathfrak{U}_{[M]}:=\sqq^{-\dim\End_{\bfk Q}(M)+\frac{1}{2}\langle M,M\rangle_Q}[M].\]
	This is the image of $\mathfrak{U}_\lambda$ under the specialization map $v\mapsto\sqq$, where $\lambda\in\mathfrak{P}$ is such that $M_q(\lambda)=M$. We can also define the $\diamond$-action of $\widetilde{\ct}(\bfk Q,\varrho)$ on $\widetilde{\ch}(\bfk Q,\varrho)$ by 
	\[
	\K_\alpha\diamond[M]=v^{\frac{1}{2}(\alpha-\varrho\alpha,\,\dimv M)_Q}\K_\alpha\ast[M].
	\]
	
	Let $\ct=\{X\in \mod(\Lambda^\imath) \mid \Hom(X,S_\ell\oplus S_{\btau \ell})=0\}$. Then any $\bfk Q$-module $M$ is of the form 
	\begin{align}
		\label{eq:M-form}
		M\cong \begin{cases}
			X_M\oplus S_\ell^{\oplus a}\oplus S_{\btau \ell}^{\oplus b}& \text{ if }\varrho\ell\neq \ell,
			\\
			X_M\oplus S_\ell^{\oplus a} &\text{ if }\btau \ell=\ell,
		\end{cases}
	\end{align}
	where $X_M\in \ct\cap \mod(\bfk Q)$.
	
	\begin{proposition}
		\label{Gamma_l on H_lambda}
		Let $M\in\mod(\bfk Q)$ be of the form \eqref{eq:M-form}, and define
		\[N=\begin{cases}
			F_{\ell}^+(X_M)\oplus (S'_{\varrho \ell})^{\oplus a}\oplus (S'_{ \ell})^{\oplus b}& \text{if $\varrho \ell\neq\ell$},\\
			F_\ell^{+}(X_M)\oplus (S_{\ell}')^{\oplus a}& \text{if $\varrho \ell=\ell$}.
		\end{cases}\]
		Then we have 
		\begin{align}
			\Gamma_\ell(\mathfrak{U}_{[M]})=\begin{cases}
				[\K_{a\alpha_\ell+b\alpha_{\varrho\ell}}]^{-1}\diamond \mathfrak{U}_{[N]} &\text{ if }\varrho \ell\neq \ell,
				\\
				[\K_{a\alpha_\ell}]^{-1}\diamond \mathfrak{U}_{[N]}&\text{ if }\varrho\ell=\ell.
			\end{cases}
		\end{align}
	\end{proposition}

	The result of Proposition~\ref{Gamma_l on H_lambda} can be made generically, and this implies $\Gamma_\ell$ preserves the rescaled Hall bases, that is, it maps the basis $\{\K_\alpha\diamond \mathfrak{U}_\lambda\mid\alpha\in\Z^\I,\lambda\in\mathfrak{P}\}$ for $\widetilde{\ch}(Q,\varrho)$ to the one for $ \widetilde{\ch}(r_\ell Q,\varrho)$.
	
	\begin{corollary}\label{iHA reflection functor of dCB}
		The isomorphism $\Gamma_\ell:\widetilde{\ch}(Q,\varrho)\rightarrow \widetilde{\ch}(r_\ell Q,\varrho)$ maps the dual canonical basis of $\widetilde{\ch}(Q,\varrho)$ to the dual canonical basis of $\widetilde{\ch}(r_\ell Q,\varrho)$.
	\end{corollary}

	The Fourier transform of $\imath$Hall algebras is constructed in \cite[\S5]{LP25}.
	Let $(Q,\btau)$ be a Dynkin $\imath$quiver, and $(Q',\btau)$ be the Dynkin $\imath$quiver constructed from $(Q,\btau)$ by reversing the arrows in a subset $\ce$ of $Q_1$ satisfying $\btau(\ce)=\ce$. Then the Fourier transform $\Phi_{Q',Q}:\widehat{\ch}(\bfk Q,\btau)\rightarrow \widehat{\ch}(\bfk Q',\btau)$ is a canonical isomorphism of $\Q(v^{1/2})$-algebras, which sends
	$[S_i]\mapsto [S_i']$, $[K_i]\mapsto [K_i']$, for $i\in\I$.
	
	\begin{corollary}[cf. {\cite[Theorem 6.4]{LP26a}}]
		The isomorphism $\Phi_{Q',Q}:\widehat{\ch}(\bfk Q,\btau)\rightarrow \widehat{\ch}(\bfk Q',\btau)$ preserves their dual canonical bases.
	\end{corollary}
	
	\begin{definition}[Dual canonical basis of $\imath$quantum groups]
		\label{def:dual CB}
		The dual canonical basis for $\widetilde{\ch}(Q,\varrho)$ is transferred to a basis for $\tUi$ via the isomorphism in Lemma~\ref{lem:Hall-iQG}, which are called the {\em dual canonical basis} for $\tUi$.
	\end{definition}
	
	The dual canonical basis of $\tUi$ can be equivalently defined by using the basis \eqref{eq:bases dual} of $\tRi$, and in particular it is integral and positive. Similarly, we can define the dual canonical basis for $\hUi$. The following theorem follows directly from Corolary~\ref{iHA reflection functor of dCB}. 
	
	\begin{theorem}\label{iQG dCB braid invariant}
		The dual canonical basis of $\tUi$ is invariant (as a set) under braid group actions.
	\end{theorem}
	
	\subsection{Dual canonical bases of quantum groups}
	\label{subsec:dCB-QG}
	
	In \cite{BG17} Berenstein and Greenstein defined a basis of $\hU$ from the dual canonical bases of $\U^+$ and $\U^-$, called the double canonical basis. In this section we recall their construction and prove that this basis coincides with the dual canonical basis of $\hU$. 
	
	
	First, following \cite{BG17} we define the quantum Heisenberg algebras $\mathcal{H}^{\pm}$ by 
	\[\mathcal{H}^+=\hU/\langle K_i'\mid i\in\I\rangle,\quad \mathcal{H}^-=\hU/\langle K_i\mid i\in\I\rangle.\]
	Let $\mathbf{K}^+$ (resp. $\mathbf{K}^-$) be the submonoid of $\hU$ generated by the $K_i$ (respectively, the $K_i'$), $i\in\I$. Then we have triangular decompositions
	\[\mathcal{H}^+=\mathbf{K}^+\otimes\hU^-\otimes \hU^+,\quad \mathcal{H}^-=\mathbf{K}^-\otimes\hU^-\otimes \hU^+.\]
	which induce natural embeddings of vector spaces
	\begin{align*}
		\iota_+&:\mathcal{H}^+=\mathbf{K}^+\otimes\hU^-\otimes \hU^+\hookrightarrow \hU=\mathbf{K}^-\otimes(\mathbf{K}^+\otimes\hU^-\otimes \hU^+),\\
		\iota_-&:\mathcal{H}^-=\mathbf{K}^-\otimes\hU^+\otimes \hU^-\hookrightarrow \hU=\mathbf{K}^+\otimes(\mathbf{K}^-\otimes\hU^+\otimes \hU^-)
	\end{align*}
	spliting the canonical projections $\hU\rightarrow\mathcal{H}^+$ and $\hU\rightarrow\mathcal{H}^-$.
	
	Let $\widetilde{\mathbf{B}}^{\pm}$ be the rescaled dual canonical basis of $\hU^{\pm}$ defined in \S\ref{dCB of HA subsec}. Recall from Lemma~\ref{QG bar-involution def} that there is a bar-involution on $\hU$ defined by $\ov{v}=v^{-1}$, $\ov{E}_i=E_i$, $\ov{F}_i=F_i$ and $\ov{K}_i=K_i$, $\ov{K}'_i=K_i'$ for $i\in\I$.
	
	Let $\Gamma=\N^\I\times\N^\I$ and set $\alpha_{+i}=(\alpha_i,0)$, $\alpha_{-i}=(0,\alpha_i)$. We define a new degree function on $\hU$ by setting
	\[\deg_\Gamma(E_i)=\alpha_{+i},\quad \deg_\Gamma(F_i)=\alpha_{-i},\quad \deg_\Gamma(K_i)=\deg_\Gamma(K_i')=\alpha_{+i}+\alpha_{-i}.\]
	It is easily seen that $\hU$ becomes a $\Gamma$-graded algebra. Using this degree function we define an action $\diamond$ of the algebra $\hU^0$ on $\hU$ via
	\begin{align*}
		K_i\diamond x=v^{-\frac{1}{2}\check{\alpha}_i(\deg_\Gamma(x))}K_ix,\quad K_i'\diamond x=v^{\frac{1}{2}\check{\alpha}_i(\deg_\Gamma(x))}K_i'x
	\end{align*}
	where $\check{\alpha}_i\in\Hom_\Z(\Gamma,\Z)$ is defined by $\check{\alpha_i}(\alpha_{\pm i})=\pm c_{ij}$ and $x$ is homogeneous. This action is characterized by the following property: 
	\begin{equation}\label{QG diamond action char}
		\ov{K\diamond x}=K \diamond \ov{x},\quad  K\in\hU^0,x\in\hU.
	\end{equation}
	Note that the $\diamond$-action as well as the bar-involution factors through to a $\mathbf{K}^{\pm}$-action and an anti-involution on $\mathcal{H}^{\pm}$ via the canonical projection $\hU\rightarrow\mathcal{H}^{\pm}$, and (\ref{QG diamond action char}) still holds. 
	
	
	\begin{theorem}[\text{\cite[Main Theorem 1.3]{BG17}}]
		\label{thm:doubleCB-H+}
		For any $(b_-,b_+)\in\widetilde{\mathbf{B}}^-\times\widetilde{\mathbf{B}}^+$, there is a unique element $b_-\circ b_+\in\mathcal{H}^+$ fixed by $\bar{\cdot}$ and satisfying
		\[b_-\circ b_+-b_-b_+\in\sum v\Z[v]K\diamond(b'_-b'_+)\]
		where the sum is taken over $K\in\mathbf{K}^+\setminus\{1\}$ and $b'_{\pm}\in\widetilde{\mathbf{B}}^{\pm}$ such that $\deg_\Gamma(b_-b_+)=\deg_\Gamma(K)+\deg_\Gamma(b'_-b'_+)$. The basis $\{K\diamond(b_-\circ b_+)\mid K\in\mathbf{K}^+,b_{\pm}\in\widetilde{\mathbf{B}}^{\pm}\}$ is called the double canonical basis of $\mathcal{H}^+$.
	\end{theorem} 
	
	\begin{theorem}[\text{\cite[Main Theorem 1.5]{BG17}}]
		\label{thm:doubleCB-U +}
		For any $(b_-,b_+)\in\widetilde{\mathbf{B}}^-\times\widetilde{\mathbf{B}}^+$, there is a unique element $b_-\bullet b_+\in\hU$ fixed by $\bar{\cdot}$ and satisfying
		\[b_-\bullet b_+-\iota_+(b_-\circ b_+)\in \sum v^{-1}\Z[v^{-1}]K\diamond\iota_+(b'_-\circ b'_+)\]
		where the sum is taken over $K\in\hU^0\setminus\mathbf{K}^+$ and $b'_{\pm}\in\widetilde{\mathbf{B}}^{\pm}$ such that $\deg_\Gamma(b_-b_+)=\deg_\Gamma(K)+\deg_\Gamma(b'_-b'_+)$. The basis $\{K\diamond(b_-\bullet b_+)\mid K\in\hU^0,b_{\pm}\in\widetilde{\mathbf{B}}^{\pm}\}$ is called the double canonical basis of $\hU$.
	\end{theorem}

	It is proved in \cite{LP26b} that the double canonical basis of $\tU$ actually coincides with the dual canonical basis. More precisely, we have the following theorem.
	
	\begin{theorem}[{\cite[Theorem 4.16]{LP26b}}]
		\label{iCB is DCB}
		The double canonical basis for $\hU$ coincides with the dual canonical basis.
	\end{theorem}

	Using Theorem~\ref{iCB is DCB}, we can answer many conjectures proposed by \cite{BG17}. For example, the following results are already proved for the dual canonical basis of $\tUi$, so they are also valid for the double canonical basis; see \cite[\S4]{LP26b}.
	
	\begin{corollary}[{\cite[Conjecture 1.15]{BG17}}; see Theorem \ref{iQG dCB braid invariant}]\label{coro: QG double CB braid}
		The double canonical basis of $\tU$ is invariant under braid group actions.
	\end{corollary}
	
	\begin{corollary}[{\cite[Conjecture 1.21]{BG17}}
		]\label{coro: QG double CB positive}
		For any $(b_-,b_+)\in\widetilde{\mathbf{B}}^-\times\widetilde{\mathbf{B}}^+$, the transition coefficients of $b_-b_+$ (or $b_+b_-$) with respect to the double canonical basis of $\hU$ belong to $\N[v,v^{-1}]$.
	\end{corollary}
	
	\begin{corollary}
		The structure constants of double canonical basis are in $\N[v^{\frac{1}{2}},v^{-\frac{1}{2}}]$. 
	\end{corollary}

	In \cite{BG17} the authors defined $\Q(v^{\frac{1}{2}})$-linear anti-involutions ${\cdot}^*:\hU\rightarrow\hU$ and ${\cdot}^t:\hU\rightarrow\hU$ via
	\begin{gather*}
		E_i^*=E_i,\quad F_i^*=F_i,\quad K_i^*=K_i',\quad (K_i')^*=K_i. 
	\end{gather*}
	It is conjectured in \cite{BG17} that ${\cdot}^*$ also preserves the double canonical basis. Using Theorem~\ref{iCB is DCB} we can now verify this conjecture.

	\begin{proposition}[{\cite[Conjecture 1.11]{BG17}}]
		\label{QG dCB invariant under *}
		The double canonical basis of $\hU$ is invariant under ${\cdot}^*$. More precisely, we have
		\[(K\diamond(b_-\bullet b_+))^*=K^*\diamond (b_-^*\bullet b_+^*).\]
	\end{proposition}

	\section{Examples: rank $1$}
	\label{sec:example}
	
	\subsection{Quantum group of rank $1$}

	In this subsection we determine the dual canonical basis of $\tU_v(\mathfrak{sl}_2)$, the algebra generated by $E,F,K,K'$. For this we need to consider the double of the quiver with one vertex $1$ and no edges. The quiver of the regular NKS category $\mcr$ is 
	\[
	\xymatrix{\tS_1\ar[d]^{\alpha_1} & \sigma(\tS_1) \ar[l]_{\beta_1} \\
		\sigma(\Sigma \tS_1) \ar[r]^{\beta_2}& \Sigma \tS_1\ar[u]^{\alpha_2} }
	\]
	subject to $\beta_2\alpha_1=0,\beta_1\alpha_2=0$. The quiver of the singular NKS category $\cs$ is 
	$$\xymatrix{\sigma(\tS) \ar@<0.5ex>[r]^{\alpha} & \sigma(\Sigma\tS) \ar@<0.5ex>[l]^{\beta}}$$ with $\beta\alpha=0=\alpha\beta$. 
	We have 
	\[\N^{\mathcal{R}_0-\mathcal{S}_0}=\N\e_{\tS_1}+\N\e_{\Sigma\tS_1},\quad \N^{\mathcal{S}_0}=\N\e_{\sigma\tS_1}+\N\e_{\sigma\Sigma\tS_1}.\]
	For vectors $\bv\in\N^{\mathcal{R}_0-\mathcal{S}_0}$ and $\bw\in\N^{\mathcal{S}_0}$, we set
	\[\bv_1=\bv(\tS_1),\quad \bv_2=\bv(\Sigma\tS_1),\quad \bw_{1}=\bw(\sigma\tS_1),\quad \bw_2=\bw(\sigma\Sigma\tS_1).\]
	
	Then the following multiplication formula is proved in \cite{LP26b}:
	
	\begin{proposition}
		For $a,b\geq 0$, denote by $c=\max(a,b)$, $d=\min(a,b)$. Then 
		\begin{align*}
			E^aF^b=
			\sum\limits_{0\leq \bv_1+\bv_2\leq d}v^{c(\bv_2-\bv_1)}\qbinom{d+1}{\bv_1}\qbinom{d+1}{\bv_2}\frac{[d+1-\bv_1-\bv_2]}{[d+1]}L(\bv,a\e_{\sigma\tS_1}+b\e_{\sigma\Sigma\tS_1}). 
		\end{align*}
	\end{proposition}
	
	Using this one can solve for the explicit expression of $L(\bv,\bw)$:
	
	\begin{proposition}\label{prop:rank I L(vw) formula}
		For any strongly $l$-dominant pair $(\bv,\bw)$, we have
		\begin{align*}
			L(\bv,\bw)=\sum_{k=l}^{n}(-1)^{k-l} \sum_{\substack{\bv'_1\geq\bv_1,\bv'_2\geq\bv_2 \\ \bv'_1+\bv'_2=k}} 
			&v^{f_{\bw}(\bv,\bv')} \qbinom{n-\bv'_2-\bv_1}{\bv'_1-\bv_1} \qbinom{n-\bv'_1-\bv_2}{\bv'_2-\bv_2}
			\\
			&\times E^{\bw_1-k}F^{\bw_2-k}K^{\bv'_1}K'^{\bv'_2},
		\end{align*}
		where $n=\min\{\bw_1,\bw_2\}$, $l=\bv_1+\bv_2$, and
		\[f_{\bw}(\bv,\bv')=(\bw_1-\bw_2)(\bv_1-\bv_2)+(n+1-\bv'_1-\bv'_2)[(\bv_1-\bv_2)-(\bv'_1-\bv'_2)].\]
	\end{proposition}
	
	From Theorem \ref{iCB is DCB} we see that the dual canonical basis of $\tU(\mathfrak{sl}_2)$ coincides with the double canonical basis defined in \cite{BG17}. In particular, by \cite[(1.1)]{BG17}, $\widetilde{\mathbf{B}}$ equals to 
	\begin{align*}
		\{v^{(a_+-a_-)(m_--m_+)}K^{a_+}K'^{a_-} F^{m_-}C^{(m_0)}E^{m_+}\mid a_{\pm}\in\Z,m_{\pm},m_0\in\N,\min(m_+,m_-)=0\},
	\end{align*}
	where $C^{(0)}=1$, $C^{(1)}=C=FE-vK-v^{-1}K'$ and 
	\begin{align}
		\label{def:Cn}
		C^{(m+1)}=CC^{(m)}-KK'C^{(m-1)},\qquad m\geq1.
	\end{align}
	Note that $C$ is central in $\tU$ and elements in $\widetilde{\mathbf{B}}$ are bar invariant, i.e., $\ov{u}=u$ for any $u\in  \widetilde{\mathbf{B}}$. Moreover, $\widetilde{\mathbf{B}}$ contains $\widetilde{\mathbf{B}}^{\pm}$.

	\subsection{$\imath$Quantum group of rank $1$}
	
	In this subsection we give the dual canonical basis of $\tUi_v(\mathfrak{sl}_2)$, the commutative algebra generated by $B,\K$, where $\K=v\tk=vKK'$. The regular NKS category $\mcr^\imath$ is
	\[
	\begin{tikzcd}
		\tS_1\ar[r,shift left=2pt,"\alpha"]& \sigma(\tS_1)\ar[l,shift left=2pt,"\beta"]
	\end{tikzcd}
	\]
	subject to $\alpha\beta=0$. The singular NKS category $\cs^\imath$ is isomorphic to $k[X]/( X^2)$.
	
	We will identify $\tUi(\mathfrak{sl}_2)$ with $\tRi$ in the following, in particular, the generators $B=L(0,\e_{\sigma\tS_1})$ and $\K=L(\bv^1,\bw^1)$. 
	
	By definition, we have
	\[\N^{\mathcal{R}^\imath_0-\mathcal{S}^\imath_0}=\N\e_{\tS_1},\quad \N^{\mathcal{S}^\imath_0}=\N\e_{\sigma \tS_1}.\] 
	So we may write $(k,m)$ for $l$-dominant pairs in $\mathcal{R}^\imath$.
	
	By \cite[\S7]{LP26a}, the \emph{dual $\imath$canonical basis} of $\tUi$ is 
	\[L(k,m)=\sum_{j=k}^{\lfloor\frac{m}{2}\rfloor}(-1)^{j-k}\binom{m-k-j}{m-2j}B^{m-2j}\K^j.\]
	Conversely, we have
	\[B^a\K^b=\sum_{i=0}^{\lfloor\frac{a}{2}\rfloor}C_{i,a-1}L(i+b,a+2b),\]
	where $C_{i,a-1}=\binom{a-1}{i}-\binom{a-1}{i-2}$.

	\section{Discussions and open problems}
	\label{sec:open}

	\subsection{Dual canonical bases under embedding}
	Recall that $\tUi$ is defined to be a subalgebra of $\tU$ and we have an embedding $\iota:\tUi\to\tU$. We have defined the dual canonical basis of $\tUi$, which we may denote by $\mathcal{B}^\imath$. Since $\tU$ can be regarded as an $\imath$quantum group of diagnoal type (see Example~\ref{ex:QGvsiQG}), its dual canonical basis is also defined, and we denote it by $\mathcal{B}$. It is then natural to consider the image of $\mathcal{B}^\imath$ under $\iota$, and we conjectural that transition coefficients from the dual canonical basis of $\tUi$ to that of $\tU$ belong to $\N[v,v^{-1}]$:
	
	\begin{conjecture}
		For each $b^\imath\in\mathcal{B}^\imath$, we have 
		\[\iota(b^\imath)=\sum_{b\in\mathcal{B}}a_{b,b^\imath}b\]
		with $a_{b,b^\imath}\in\N[v,v^{-1}]$.
	\end{conjecture}
	
	This conjecture has already been proved for $\tUi_v(\mathfrak{sl}_2)$, see \cite{CZ25}.

	\subsection{Dual canonical bases under coproduct}
	
	The dual canonical basis $\mathcal{B}$ of $\tU$ is defined via the isomorphism $\tilde{\kappa}:\tU\stackrel{\simeq}{\rightarrow}\tR$, and since the multiplication of $\tR$ is defined using the restriction functor $\Delta^{\bw}_{\bw^1,\bw^2}$, we know that the structural constants of $\mathcal{B}$ for multiplication live in $\N[v,v^{-1}]$. However, $\tU$ is a bialgebra and we can also consider the structural constants of $\mathcal{B}$ for coproduct. We conjecture that the structural constants for coproduct also belong to $\N[v,v^{-1}]$:
	
	\begin{conjecture}
		For each $b\in\mathcal{B}$, we have 
		\[\Delta(b)=\sum_{b_1,b_2\in\mathcal{B}}c^b_{b_1,b_2}b_1\otimes b_2\]
		where $c^b_{b_1,b_2}\in \N[v,v^{-1}]$.
	\end{conjecture}

	We know $\tUi$ is a right coideal subalgebra of $\tU$, that is, $\Delta(\tUi)\subset \tUi\otimes\tU$. Then $\mathcal{B}^\imath\otimes\mathcal{B}=\{b^\imath\otimes b\mid b^\imath\in\mathcal{B}^\imath,b\in\mathcal{B}\}$ is a basis of $\tUi\otimes \tU$. 
	
	\begin{conjecture}
		For each $b\in\mathcal{B}$, we have 
		\[\Delta(\iota(b^\imath))=\sum_{b_1,b_2\in\mathcal{B}}c^{b^\imath}_{b_1^\imath,b_2}b_1^\imath\otimes b_2\]
		where $c^{b^\imath}_{b_1^\imath,b_2}\in \N[v,v^{-1}]$.
	\end{conjecture}

	Again, these two conjectures are proved for $\tU_v(\mathfrak{sl}_2)$; see \cite{CZ25}.
	\subsection{ For $\imath$quantum group of type ${\rm AIII}_{2r}$}
	
	We have constructed $\imath$Hall algebra and dual canonical basis from an $\imath$quiver $(Q,\varrho)$. It turns out that not every Dynkin diagram with involution $\varrho$ can be equipped with an orientation $\Omega$ such that $\varrho$ becomes an involution of quiver. The only exception for type ADE is the so-called type ${\rm AIII}_{2r}$ case, whose Dynkin diagram and involution $\varrho$ are given as follows:
	\begin{center}\setlength{\unitlength}{0.6mm}
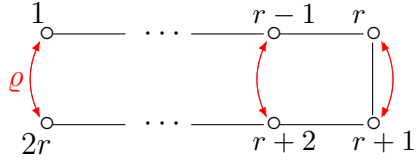
		\begin{picture}(70,35)(0,5)
			
			
			\put(0,10){$\circ$}
			\put(0,30){$\circ$}
			\put(50,10){$\circ$}
			\put(50,30){$\circ$}
			\put(72,10){$\circ$}
			\put(72,30){$\circ$}
			\put(-4,5){$2r$}
			\put(-2,34){${1}$}
			\put(47,6){\small $r+2$}
			\put(47,34){\small $r-1$}
			\put(69,6){\small $r+1$}
			\put(69,34){\small $r$}
			
			\put(3,11.5){\line(1,0){16}}
			\put(3,31.5){\line(1,0){16}}
			\put(23,10){$\cdots$}
			\put(23,30){$\cdots$}
			\put(33.5,11.5){\line(1,0){16}}
			\put(33.5,31.5){\line(1,0){16}}
			\put(53,11.5){\line(1,0){18.5}}
			\put(53,31.5){\line(1,0){18.5}}
			
			\put(73.5,13.6){\line(0,1){16}}
			
			\color{red}
			\put(-7,20){$\varrho$}
			\qbezier(0,13.5)(-4,21.5)(0,29.5)
			\put(-0.25,14){\vector(1,-2){0.5}}
			\put(-0.25,29){\vector(1,2){0.5}}
			
			\qbezier(50,13.5)(46,21.5)(50,29.5)
			\put(49.75,14){\vector(1,-2){0.5}}
			\put(49.75,29){\vector(1,2){0.5}}

			\qbezier(76,13.5)(80,21.5)(76,29.5)
			\put(75.85,14){\vector(-1,-2){0.5}}
			\put(75.85,29){\vector(-1,2){0.5}}
		\end{picture}
	\end{center}
	
	In \cite{CLPRW25}, the authors develop an $\imath$Hopf algebra approach to realize $\imath$quantum groups; compare the star product constructed in \cite{KY20,KY21}. 
	Using this construction, one can still define the dual canonical basis in this type, and Theorem~\ref{iQG dCB braid invariant} remains valid; see \cite{CLPRW25}. We conjecture that the positivity results also hold in this case:
	
	\begin{conjecture}
		\label{conj:positive}
		The positivity results for dual canonical basis still hold for type ${\rm AIII}_{2r}$. For example, the structural constants for multiplication belong to $\N[v,v^{-1}]$.
	\end{conjecture}

	\subsection{Non-simply-laced type}
	
	The $\imath$Hopf algebra approach mentioned above works for all $\imath$quantum groups of quasi-split type. In this way, we can also define the dual canonical basis for $\tUi$ of finite type, and prove that it is invariant under braid group actions; see \cite{CLPRW25}. Similarly to the canonical bases in \cite{Lus93} for non-simply-laced type, we can not expect the positivity stated in Conjecture \ref{conj:positive} to hold in general. Using the language of quivers with automorphisms, Lusztig \cite{Lus93} also gave a geometric realiation of quantum groups of non-simply-laced type. 
	
	We shall formulate valued $\imath$quiver algebras, and then use the $\imath$Hall algebras of valued $\imath$quiver algebras to realize these $\imath$quantum groups.	
	
	\begin{problem}
		Develop the construction of NKS quiver varieties further to realize $\tUi$ of non-simply-laced type. The perverse sheaves correspond to the dual canonical basis. The transition matrix coefficients of the $\imath$Hall basis with respect to the dual canonical basis is positive; cf. \cite[Conjecture 1.21]{BG17} for quantum groups. 
	\end{problem}

	\subsection{Beyond quasi-split type}

	A Satake diagram consists of a pair $(\I =\I_\bullet \cup \I_\circ, \varrho)$, where $\I =\I_\bullet \cup \I_\circ$ is a bi-colored partition of the Dynkin diagram $\I$, and $\varrho$ is a diagram involution (we allow $\varrho=\Id$), subject to suitable conditions. In finite type, the Satake diagrams are classified by Araki; cf. \cite[Table 4]{BW18b}; they are in bijection with real forms of complex simple Lie algebras.

	\begin{problem}
		Provide $\imath$Hall algebra realization and NKS geometric realization of $\tUi$ associated to a general Satake diagram of finite type, and construct the dual canonical basis for $\tUi$.
	\end{problem}
	
	For general Satake diagrams of type ADE, we expect the dual canonical bases also have the same integral and positive properties as the quasi-split type.



\end{document}